\documentclass[aip,jcp,reprint,floatfix,a4paper]{revtex4-1}  
\usepackage{cmap}
\usepackage[T1]{fontenc}
\usepackage{times}
\usepackage{helvet}
\usepackage{fixmath}
\DeclareSymbolFont{UPM}{U}{eur}{m}{n}
\DeclareMathSymbol{\partial}{0}{UPM}{"40}
\usepackage{graphicx}
\usepackage{amsmath}
\usepackage{amssymb}
\usepackage{geometry}
\usepackage{hyperref}
\hypersetup{
        colorlinks=true,
        citecolor=black,
        filecolor=black,
        linkcolor=black,
        urlcolor=black
}

\usepackage[dvipsnames]{xcolor}
\usepackage[stretch=15,shrink=15,step=1]{microtype}
\SetProtrusion{encoding={*},family={*},series={*},size={6,7}}
              {1={ ,550},2={ ,300},3={ ,300},4={ ,300},5={ ,300},
               6={ ,300},7={ ,400},8={ ,300},9={ ,300},0={ ,300}}

\usepackage[nodayofweek]{datetime}
\newdateformat{myDate}{\THEDAY\ \monthname[\THEMONTH] \THEYEAR}

\usepackage{bm}
\usepackage{psfrag}
\usepackage{placeins}
\usepackage{csquotes}
\usepackage{mathtools}
\usepackage{color}
\usepackage{bbm}
\makeatletter \newcommand \listoftodos{\section*{Todo list} \@starttoc{tdo}}
\newcommand\l@todo[2]
 {\par\noindent \textbf{#2}:  #1\par\vspace{0.2cm}} \makeatother

\setcitestyle{authoryear,round}

\begin{document}

\title[]{3D cut-cell modelling for high-resolution atmospheric simulations}

\author{H. Yamazaki}
\email[Correspondence author. E-mail: ]{h.yamazaki@imperial.ac.uk}
\affiliation{
Department of Mathematics, Imperial College London, London, UK
}%

\author{T. Satomura}%
\affiliation{ 
Division of Earth and Planetary Sciences, Graduate School of Science, Kyoto University, Kyoto, Japan
}%

\author{N. Nikiforakis}
\affiliation{%
Department of Physics, Cavendish Laboratory, University of Cambridge, Cambridge, UK
}%

\date{\today}

\begin{abstract}
Owing to the recent, rapid development of computer technology, the resolution of atmospheric numerical models has increased substantially. With the use of next-generation supercomputers, atmospheric simulations using horizontal grid intervals of $O$(100) m or less will gain popularity. At such high resolution more of the steep gradients in mountainous terrain will be resolved, which may result in large truncation errors in those models using terrain-following coordinates. 
In this study, a new 3D Cartesian coordinate non-hydrostatic atmospheric model is developed. A cut-cell representation of topography based on finite-volume discretization is combined with a cell-merging approach, in which small cut-cells are merged with neighboring cells either vertically or horizontally. In addition, a block-structured mesh-refinement technique is introduced to achieve a variable resolution on the model grid with the finest resolution occurring close to the terrain surface.
The model successfully reproduces a flow over a 3D bell-shaped hill that shows a good agreement with the flow predicted by the linear theory. The ability of the model to simulate flows over steep terrain is demonstrated using a hemisphere-shaped hill where the maximum slope angle is resolved at 71$^{\circ}$. The advantage of a locally refined grid around a 3D hill, with cut-cells at the terrain surface, is also demonstrated using the hemisphere-shaped hill. The model reproduces smooth mountain waves propagating over varying grid resolution without introducing large errors associated with the change of mesh resolution. At the same time, the model shows a good scalability on a locally refined grid with the use of OpenMP.
\end{abstract}


\maketitle

\section{Introduction}\label{s:intro}
One of the pressing concerns of next-generation high-resolution atmospheric modeling is the accurate treatment of terrain. The continuous increase in computer power and the associated increase in model resolution has resulted in the resolution of steeper and more complex features in the terrain. These variations in the bottom surface of the atmosphere not only have a significant influence on the local dynamics near the surface but can also affect the global circulations \citep{McFarlane1987}. Although the resolutions of 1--20 km are used in today's operational models, higher-resolution simulations at horizontal grid intervals of $O(100)$ m or less will gain popularity with the use of next-generation supercomputers \citep[e.g.,][]{MiyamotoEA2013}. Therefore it becomes more important for the future high-resolution models to implement a robust method of representing topography for steep gradients and complex geometries. 

For many years, the common choice for the representation of topography in atmospheric models has been the terrain-following vertical coordinates based either on pressure \citep[e.g.,][]{Philips1957, Simmons+Burridge1981}, or on height \citep[e.g.,][]{Gal-Chen+Somerville1975}. In the terrain-following coordinates, the vertical model levels follow the shape of the terrain at the bottom and gradually revert to horizontal surfaces with increasing height above the surface. The main advantage of this approach is that the imposition of the lower boundary condition is straightforward for arbitrary topography. In addition, the terrain-following coordinates are suitable for coupling with boundary layer parameterizations because a high near-ground resolution is easily achieved by increasing the number of model levels near the bottom boundary. Although the terrain-following coordinates have proven effective for a wide range of applications, large truncation errors may arise in computing the horizontal gradients, particularly in the presence of steep terrain \citep{ThompsonEA1985, Satomura1989}. It is recognized that the most critical error lies in the discretization of the horizontal pressure gradient term, which can induce spurious circulations over mountainous topography \citep{Janjic1989} and even numerical instability if the mountains are steep enough \citep{Zangl2012}. 

There are many ongoing efforts to alleviate the disadvantages of the terrain-following coordinates. A substantial improvement has been made by the specification of the vertical coordinate which removes the influence of the topography as fast as possible with height \citep{ScharEA2002, LeuenbergerEA2010, Klemp2011}. This successfully reduced the errors at upper levels where the coordinates are much smoother than classical hybrid coordinates. Another improvement lies in a better treatment of the horizontal gradient of pressure \citep{Klemp2011, Zangl2012} as well as of diffusion \citep{Zangl2002}. While these studies have improved the accuracy of the terrain-following coordinates substantially, a conclusion has not been reached on whether the terrain-following models would be accurate enough for future generation high-resolution models.

This study examines methods based on the use of Cartesian coordinates as an alternative means of representing topography. In these methods, the terrain is directly incorporated into a regular rectangular grid without using a coordinate transformation. Hereafter in this study, a regular rectangular grid system is referred to as a `Cartesian grid' and a representation method of topography based on a Cartesian grid is referred to as a `Cartesian-grid method'. Since the model levels are kept horizontal throughout the domain, Cartesian-grid methods resolve the imbalances which occur on a terrain-following grid in the discretization of the horizontal gradients. However the imposition of the lower boundary condition can be complicated in Cartesian-grid methods because the terrain surface does not normally coincide with the grid lines. 

The step-mountain method is a Cartesian-grid method with a straightforward imposition of the lower boundary condition by approximating the terrain surface as a piecewise constant function along the grid lines \citep[e.g.,][]{Bryan1969, MeisingerEA1988}. Because of the lack of accuracy in the stepwise approximation of topography, however, this method introduces serious errors at step-corners into flow patterns \citep{Gallus+Klemp2000} and therefore turned out to be ill suited for high-resolution simulations over mountains \citep{ScharEA2002, Zangl2003}. The partial-cell method provides another stepwise approximation where the heights of the steps are adjusted to those of the topography \citep[e.g.,][]{Semtner+Mintz1977, Maier-Reimer+Mikolajewicz1992}. \cite{AdcroftEA1997} showed that the partial-cell method successfully reduced the errors associated with the stepwise boundary over gently-sloping terrain. However the errors can still be large at the step-corners over steep terrain where the topographic height varies substantially in a horizontal grid length \citep{Yamazaki+Satomura2010}. 

A smoother and more precise representation of the terrain is achieved by allowing linear variation of the boundary within a cell, resulting in various shapes of cells that are cut by the terrain surface. Finite-volume discretization of the governing equations assures conservation of model variables on those irregularly-shaped cells as well as on regular uncut cells. This finite-volume based Cartesian-grid method is refered to as the cut-cell (or shaved-cell) method, and it is the approach explored in this paper.

The cut-cell method has been most popular in the field of computational fluid dynamics for simulating flows with complex geometry \citep[e.g.,][]{Quirk1994, PemberEA1995, UdaykumarEA1996, YeEA1999}. After initial implementation in an ocean model \citep{AdcroftEA1997}, the cut-cell method has been examined for accuracy and robustness in atmospheric models over the past dozen years \citep[e.g.,][]{SteppelerEA2002, SteppelerEA2006, Walko+Avissar2008, Yamazaki+Satomura2008, Yamazaki+Satomura2010, KleinEA2009, Jebens2011, LockEA2012, GoodEA2013}. Many applications to well-known idealized flows in these studies suggest that the cut-cell method does not suffer from the problems reported by \cite{Gallus+Klemp2000} in the step-mountain method and can reproduce smooth mountain waves over the terrain. 

Comparisons of 2D flow results using the cut-cell method to results from terrain-following models have been made in various studies. For example, \cite{Yamazaki+Satomura2008, Yamazaki+Satomura2010} demonstrated that the cut-cell method successfully eliminated the spurious vertical velocity modes that occurred in the vicinity of steep slopes in a terrain-following model.  \cite{GoodEA2013} compared the errors in the flow aloft, where grids are fully rectangular in cut-cell models but are distorted in terrain-following models because of the influence of the underlying topography. They showed that the errors associated with a terrain-following grid are reduced when the cut-cell method is used. They also examined the robustness of the cut-cell method for steep gradients and demonstrated that it can produce stable results for flows over bell-shaped hills with aspect ratios of the height to the half-width up to 10.

The main problem associated with the use of the cut-cell method is the generation of arbitrarily small cut-cells near the boundary. Such cells lead to severe stability constraints as a result of the Courant--Friedrichs--Lewy (CFL) condition, and therefore require very small time steps. Several approaches have been introduced to resolve this small-cell problem in atmospheric cut-cell models. In a simple approach proposed by \cite{SteppelerEA2002}, the  computational volumes of cut-cells are artificially increased to those of regular uncut cells. This method was called `thin-wall' approximation because, by expanding the volume to the notionally full value but leaving the areas untouched, the terrain looks like a collection of infinitesimally thin-walls \citep{Adcroft2013}. \cite{Yamazaki+Satomura2010} use a different approach in which small cut-cells are merged with adjacent cells either vertically or horizontally. This cell-merging technique makes it possible to extend the stability limit maintaining the rigid evaluation of cut-cell volumes and areas, and thus maintaining a sharp representation of the terrain surface. In \cite{KleinEA2009}, a dimensional splitting technique was used to approximate the fluxes at cut-cell interfaces, thereby allowing the use of a full-time step defined by the regular grid. The use of implicit schemes, introduced in an atmospheric model by \cite{Jebens2011}, is another way to stabilize small cells.

Another disadvantage of the cut-cell method compared to the terrain-following approach is that to obtain a high vertical resolution at boundary layers over a wide range of topographic height can be expensive \citep{Walko+Avissar2008}. In terrain-following models, a high vertical resolution near the terrain surface is easily achieved over all topographic heights by increasing the number of model levels near the bottom boundary. However, in a Cartesian-grid based model, a substantial number of additional model levels would be required to cover all topographic heights. Around steep slopes, horizontal grid intervals must be closely spaced as well as vertical grid intervals to achieve high near-ground resolution on a Cartesian grid. A block-structured mesh refinement approach proposed in \cite{Yamazaki+Satomura2012} is one way to achieve a locally refined Cartesian grid with high computational efficiency. They demonstrated that a flow over a 2D semicircular hill was successfully reproduced on a grid locally refined around the hill with the use of cut-cells at the boundary. 

Compared to the number of 2D studies of the cut-cell representation of topography, applications of cut-cells to 3D atmospheric model studies are relatively scarce. \cite{SteppelerEA2006} investigated the impact and potential use of cut-cells in a 3D forecasting model. They demonstrated that precipitation scores and RMSE of the temperature for 1-day forecasts using a cut-cell model were improved compared to the results from the terrain-following version of the model. In the extended forecasts of 5 days, the improvements became more substantial \citep{SteppelerEA2011, SteppelerEA2013}. 

\cite{LockEA2012} focused on the capability of the cut-cell method for 3D idealized flows with a very steep slope. They showed that a potential flow over a 3D bell-shaped hill was successfully reproduced using cut-cells where the maximum slope angle of the hill was as steep as 74$^{\circ}$ from the horizontal. However no method to address the small cell problem was adopted in their cut-cell model, leaving stability and efficiency issues for practical applications.  

In this study, we propose a new 3D atmospheric model using the cut-cell method for high-resolution simulations over steep topography. To avoid the severe stability constraints from small cells, the cell-merging technique developed by \cite{Yamazaki+Satomura2010}, hereafter YS10, is extended to 3D and implemented in the model. In addition, a 3D version of the block-structured mesh refinement approach proposed in \cite{Yamazaki+Satomura2012}, hereafter YS12, is introduced to achieve high resolution near the terrain surface and also easy parallelization. Section 2 describes how the 2D methods of YS10 and YS12 are converted into 3D. Then we demonstrate the performance of the model in section 3 through numerical experiments on mountain waves. To investigate the capability of the proposed 3D cut-cell method for a wide range of slope angles, the results of flow over a 3D bell-shaped hill and a hemisphere-shaped hill are presented. In addition, the performance of the model on a locally refined grid in combination with the use of cut-cells near the boundary is examined using the hemisphere-shaped hill.

\section{Model Description}\label{s:model}

\subsection{Governing equations}\label{ss:equations}

The model solves fully compressible quasi-flux-form equations developed by \cite{Akiba2002} and \cite{Satomura+Akiba2003}, given by the following conservation equations for momentum, potential temperature and mass based on the Cartesian coordinates: 
\begin{eqnarray}
\frac{\partial\rho u}{\partial t} &=& - \nabla\cdot(\rho u\mathbf{ u}) -\frac{\partial p^{\prime}}{\partial x} + D_{u} \hspace{3pt} , \label{eq:momentum_u} \\
\frac{\partial\rho v}{\partial t} &=& - \nabla\cdot(\rho v\mathbf{ u}) -\frac{\partial p^{\prime}}{\partial y} + D_{v} \hspace{3pt} , \label{eq:momentum_v} \\
\frac{\partial\rho w}{\partial t} &=& - \nabla\cdot(\rho w\mathbf{ u}) -\frac{\partial p^{\prime}}{\partial z} - \rho^{\prime}g +D_{w} \hspace{3pt} , \label{eq:momentum_w} \\
\frac{\partial p^{\prime}}{\partial t} &=& - \frac{c_{p}R}{c_{v}p_{0}}\left(\frac{p}{p_{0}}\right)^{R/c_{p}}\left(\nabla\cdot(\rho\theta\mathbf{ u})-D_{\theta}\right) , \label{eq:heat} \\ 
\frac{\partial\rho^{\prime}}{\partial t} &=& - \nabla\cdot(\rho\mathbf{ u})  \hspace{3pt}, \label{eq:continuity}
\end{eqnarray}
where $\mathbf{u}=(u, v, w)$ is the velocity vector, $\theta$ is the potential temperature, $g$ is the acceleration due to gravity, $p$ and $\rho$ are the total pressure and total density, respectively, and the prime indicates perturbations from the hydrostatically balanced state:
\begin{eqnarray}
p &=& \bar{p}_{(z)}+p^{\prime}_{(x,y,z,t)}, \\
\rho &=& \bar{\rho}_{(z)}+\rho^{\prime}_{(x,y,z,t)}, \\
\frac{\partial\bar{p}}{\partial z} &=& -\bar{\rho}g. 
\end{eqnarray}
In Eq. (\ref{eq:heat}), $c_{p}$ and $c_{v}$ denote the specific heats at constant pressure and constant volume, respectively, $R$ is the gas constant, and $p_{0}$ is a reference pressure of $10^{5}$ Pa. The terms $D_{u}$, $D_{v}$, $D_{w}$ and $D_{\theta}$ in Eqs (\ref{eq:momentum_u})--(\ref{eq:heat}) represent source terms due to mixing and diffusion. In this study, turbulent parameterization terms based on turbulent kinetic energy \citep{Klemp+Wilhelmson1978} are used. In addition, a fourth-order artificial diffusion term is introduced in the horizontal and vertical directions to suppress the numerical noise.  Finally the system of equations is closed by the following equation of state for an ideal gas:
\begin{equation}
\theta =  \frac{p}{\rho R }\left(\frac{p_{0}}{p}\right)^{R/c_{p}}. 
\end{equation}

In the view of conservation characteristics, flux-form equations are well suited to the finite-volume discretization that is an essential part of the cut-cell method. \cite{Satomura+Akiba2003} demonstrated that the equations achieve mass conservation by simulating heat-island circulation. In addition, they designed the equations to avoid cancellation errors stemming from subtracting the hydrostatic variable ($\bar{p}$ or $\bar{\rho}$) from the nearly hydrostatic total variable ($p$ or $\rho$). This may occur in other flux-form equations \citep[e.g.,][]{SaitoEA2001, KlempEA2007}. Specifically, \cite{Satomura+Akiba2003} accomplish it by directly predicting the perturbations of the variables ($p'$ or $\rho'$).

\begin{figure*}[t]
  \centering
  \noindent\includegraphics[width=41pc,angle=0]{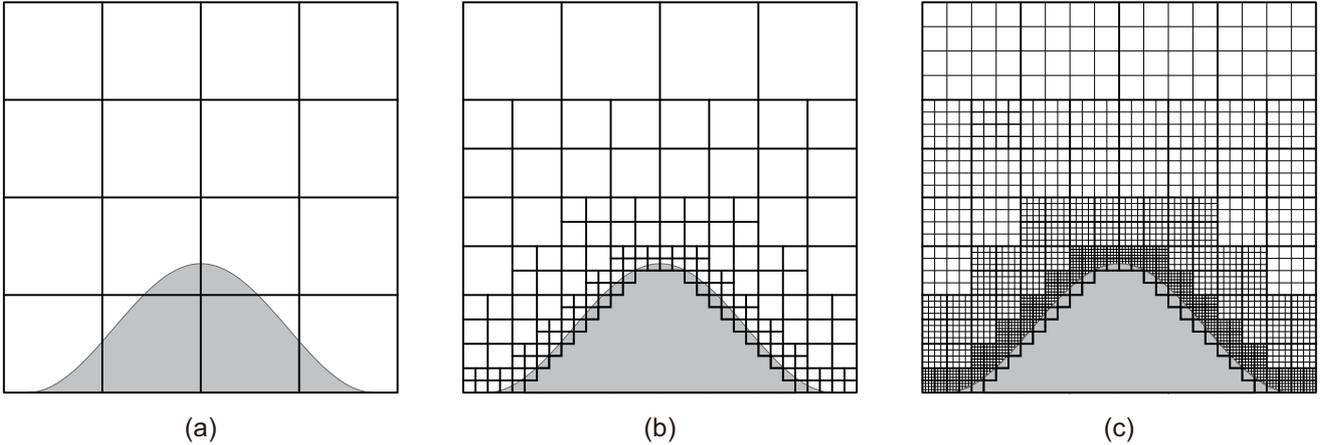}\\
  \caption{Schematics of the generation of a block-structured Cartesian mesh: (a) and (b) show the cubes before and after the refinement process, respectively, and (c) shows the cells as well as the cube boundaries. Thick and thin lines represent cube boundaries and cell boundaries, respectively. Shaded regions describe the topography.}
  \label{Fig1}
\end{figure*}

\subsection{Block-structured grid}\label{ss:blocks}

To achieve variable resolution in the Cartesian-coordinate system, the generation of the model grid begins with the creation of a block-structured grid. Here we use the Conserved Building-Cube Method (CBCM) that is originally proposed in 2D by YS12. The method is based on the Bulding-Cube Method (BCM) developed by \cite{Nakahashi2003} that has been used for the problems of computational fluid dynamics \citep[e.g.,][]{NakahashiEA2006, KimEA2007}. CBCM and BCM share the two-tiered data structure of a generated block-structured grid: tree-structured sub-domains called `cubes' and array-structured uniform mesh called `cells' in each cube.

Figure \ref{Fig1} illustrates a 2D example of the grid generation process of CBCM. First the model domain is divided into coarse equally-spaced cubes (Figure \ref{Fig1}(a)). Next, cubes that are in the vicinity of the boundary are divided into four cubes (or eight in 3D). By repeating this refinement process, the domain is divided into a number of cubes, where the size of a cube becomes small closer to the terrain surface (Figure \ref{Fig1}(b)). Here the size differences between cubes are adjusted to guarantee a uniform 2 : 1 mesh resolution at fine-coarse cube boundaries in horizontal, vertical and diagonal directions. Note that the cubes that are located completely inside of the topography are removed and are not used for the computation. Finally, a Cartesian grid of equal spacing and equal number of cells in each cube is generated (Figure \ref{Fig1}(c)). For example, each cube in Figure \ref{Fig1}(c) has $4^{2}$ cells (or $4^{3}$ in 3D) regardless of the size of the cube. Note that the generation of cut-cells in the finest cubes is described in the following subsection \ref{ss:cut-cells}.

The local grid interval is determined as
\begin{eqnarray}
	H (l) = 2^{-l}H, 
\end{eqnarray}
where the refinement level $l$ ranges from 0 (coarsest) to $l_{max}$ (finest) and the grid interval at the unrefined cube $H \equiv H (0)$. In this study, the same grid interval $H$ is used in $x$, $y$ and $z$ directions for simplicity, though it can easily be extended to different spacing in each direction. When setting up a block-structured mesh, we choose the value of $H$ so that the field far away from the topography is resolved at this resolution. Then we keep refining the cubes near the terrain surface until we obtain as high a resolution as we wish to impose close to the surface. For example, when we use $H$ = 1 km, we can get a local grid interval of 15.625 m at the refinement level 6.

A block-structured grid constructed from cubes and cells provides several attractive features for Cartesian-grid models. First, the tree-based data structure of cubes provides a locally refined grid around arbitrary topography by refining the size of the cubes near the terrain surface. The use of a uniform Cartesian mesh in each cube, at the same time, allows the direct use of any existing code based on a Cartesian grid by treating each cube as an independent computational domain. Furthermore, the same number of cells among all cubes makes the method suitable for parallel computing. Because the load balance of each cube is equivalent, high parallel efficiency can be achieved by simply distributing an equal number of cubes for each processor. Both CBCM and BCM showed good speedups with straightforward OpenMP-based parallelization \citep[][]{KimEA2007, Yamazaki+Satomura2012}.

CBCM employs a subcycling time integration that allows the use of a larger time step at coarse cubes than that used for smaller cubes. For cubes at the refinement level $l$, the time step is chosen as
\begin{eqnarray}
	\Delta t(l) = 2^{-l}\Delta t, 
\end{eqnarray}
where the global time step $\Delta t \equiv \Delta t(0)$. The equations are integrated from the cubes at the finest level. After the cells in the cubes at level $l$ are advanced in two time steps of $\Delta t(l)$, those at level $l-1$ are advanced in one time step of $\Delta t(l-1)$, followed by the information exchange between the levels. In our model, this subcycling integration is incorporated into the leap-frog time-stepping scheme and used along with the Robert-Asselin filter \citep[][]{Robert1966, Asselin1972}.

\begin{figure}[t]
  \centering
  \noindent\includegraphics[width=19pc,angle=0]{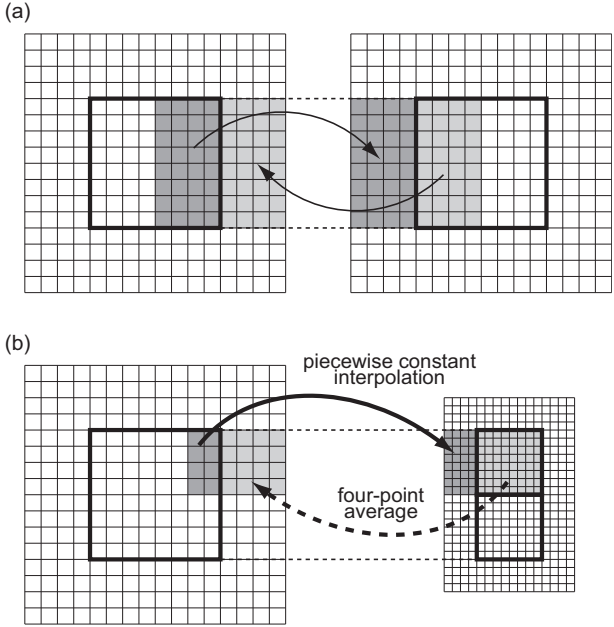}\\
  \caption{Information transfer between the adjacent cubes of (a) the same size and (b) different sizes. Thick lines represent cube boundaries, and shaded regions describe the overlapped regions of the ghost cells. Solid arrows in (a) describe information transfer between the same-size cubes. Thick solid and dotted arrows in (b) describe information transfer from a coarse cube to a fine cube, and transfer from a fine cube to a coarse cube, respectively. In this figure a cube has 8 $\times$ 8 cells with 4 ghost cells beyond the boundary.}
  \label{Fig2}
\end{figure}

Flow information is transferred between adjacent cubes through ghost cells that are added beyond the boundary of each cube. In this study, four ghost cells are added as shown in Figure \ref{Fig2}. Note that all ghost regions are at the same resolution as the inner domain of the cube. The information between the same-size cubes can be exchanged in a straightforward way because of the exact overlapping of the cells (Figure \ref{Fig2}(a)). Some interpolation methods are required at fine-coarse cube boundaries to exchange the boundary values at different resolutions (Figure \ref{Fig2}(b)). In this study, the values of ghost cells of a coarse cube are assigned by using a simple four- or eight-point average in 2D or 3D, respectively, of fine-cell values at the corresponding location in the fine cube. The ghost-cell values of a fine cube, on the other hand, are interpolated by assigning the same coarse-cell value to the corresponding four or eight ghost cells of the fine cube in 2D or 3D, respectively. Though this interpolation procedure has been demonstrated in YS12 to conserve the global second-order accuracy of the CBCM, a higher-order interpolation method can be used to improve the accuracy at fine-coarse cube boundaries \citep[e.g.,][]{JablonowskiEA2006}. 

\begin{figure}[t]
  \centering
  \noindent\includegraphics[width=19pc,angle=0]{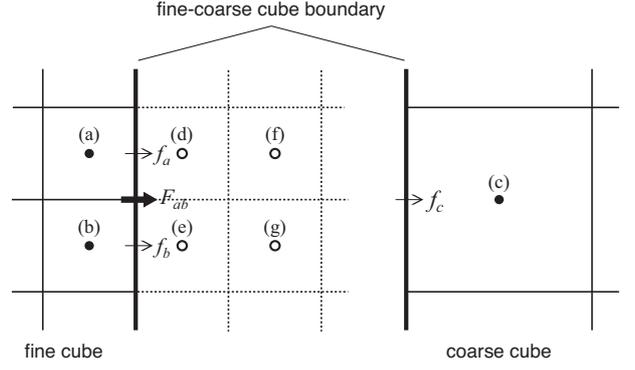}\\
  \caption{Computational cells and fluxes at a fine-coarse cube boundary. Thick and thin lines represent cube and cell boundaries, respectively. Dotted lines represent the boundaries of ghost-cells of the fine cube. Filled and open circles describe computational nodes of the cells inside the cube boundary and those of the ghost-cells, respectively. Thin and thick arrows indicate the fluxes at the cell and cube boundaries, respectively.}
  \label{Fig3}
\end{figure}

Without time interpolation, some time inconsistency may occur during the subcycling integration. For example, in a 2D case of Figure \ref{Fig3}, the coarse-cell value at (c) is assigned to the ghost cells of the fine cube (d) to (g) and can be used twice in a row without updating. This causes computational modes to the fine-coarse cube boundary due to the time suspension at the ghost region. CBCM avoids this problem by integrating equations on some of the ghost cells of the fine cube and using the updated values for the subcycling steps. The number of ghost cells on which the equations are integrated is chosen to satisfy the numerical stencils at the boundary cells inside the fine cube, such as the cells (a) and (b) in Figure \ref{Fig3}. In our model, we use a 5-point stencil in the horizontal and vertical directions to calculate fourth-order artificial diffusion terms. Therefore we integrate equations on two more cells beyond the boundary to satisfy the numerical stencils of the boundary cells with updated values. In case of Figure \ref{Fig3}, the ghost-cell values at (d) to (g) are updated as well as the values at (a) and (b) during the subcycling integration. To be able to integrate equations on two more cells beyond the boundary, four ghost cells in total are used in this study beyond each side of the cube boundary. This approach prevents the computational modes from contaminating the results inside the cubes by simply integrating equations on some extra cells beyond the boundary of each cube. On the other hand, it demands higher computational costs compared to a time interpolation scheme due to a relatively large number of ghost cells.

Another characteristic of CBCM is that it ensures global mass-conservation on a locally refined mesh. To achieve conservation with a subcycling time-stepping scheme, we must ensure that, at each fine-coarse grid interface, the numerical flux on the coarse grid equals the flux on the fine grid accumulated during the subcycling steps. CBCM facilitates this process by introducing the cube-boundary flux at fine-coarse cube boundaries. For example, in the case of Figure \ref{Fig3}, the cube-boundary flux $F_{ab}$ is defined at the same location as the coarse-cell flux $f_{c}$. During the subcycling integration of the fine cube, the fine-cell fluxes $f_{a}$ and $f_{b}$ are accumulated and stored in $F_{ab}$. Then, in the integration of the coarse cube, $F_{ab}$ overrides the coarse cell flux $f_{c}$ and is used to update the coarse-cell value at (c). In the 3D method, four fine-cell fluxes are accumulated and stored in each corresponding cube-boundary flux. This flux-matching algorithm ensures mass-conservation on the condition that the density is defined at the cell-centers, as demonstrated in YS12.

\subsection{Cut-cell configuration}\label{ss:cut-cells}

\begin{figure}[t]
  \centering
  \noindent\includegraphics[width=19pc,angle=0]{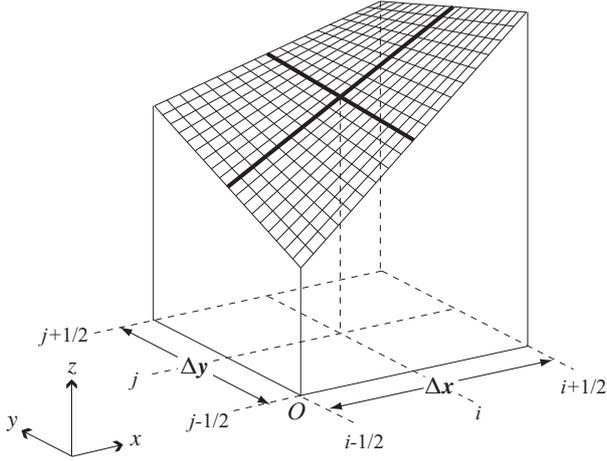}\\
  \caption{Bilinear representation of the terrain surface at the grid column ($i, j$). The corner $O$ acts as the origin of the bilinear function at the column. Thick lines indicate the boundaries at the vertical cross sections through the center of the grid column, which determine the gradients of the surface.}
  \label{Fig4}
\end{figure}

Following the generation of  a block-structured grid in the previous subsection, cut-cells are generated near the terrain surface. Because the grid is refined around the topography, the procedure of cut-cell generation is only required in the finest cubes. Following \cite{SteppelerEA2006} and \cite{LockEA2012}, the topographic boundary is represented by piecewise bilinear surfaces that are continuous at the boundaries of grid columns. First the terrain heights are specified at the four corners of the grid columns in each cube, then a unique surface for each grid column is defined by using a bilinear function of height with respect to the horizontal position ($x, y$) as
\begin{eqnarray}
	h(x,y) = m_{1}x+m_{2}xy+m_{3}y+c, \label{eq:bilin}
\end{eqnarray}
where $m_{1}$, $m_{2}$, $m_{3}$ and $c$ are constants. The function gives a linear spline at any vertical cross-sections in the $x$ or $y$ direction, as shown in Figure \ref{Fig4}. By treating the corner $O$ at ($i-1/2$, $j-1/2$) as the origin, the four heights at the corners ($i\pm1/2$, $j\pm1/2$) determine a bilinear surface at the grid column ($i$, $j$) through the coefficients  
\begin{flalign}
	m_{1_{ij}} &= \frac{h_{i+\frac{1}{2}j-\frac{1}{2}}-h_{i-\frac{1}{2}j-\frac{1}{2}}}{\Delta x},\\
	m_{2_{ij}} &= \frac{h_{i-\frac{1}{2}j-\frac{1}{2}}-h_{i-\frac{1}{2}j+\frac{1}{2}}+h_{i+\frac{1}{2}j+\frac{1}{2}}-h_{i+\frac{1}{2}j-\frac{1}{2}}}{\Delta x \, \Delta y}, \\
	m_{3_{ij}} &= \frac{h_{i-\frac{1}{2}j+\frac{1}{2}}-h_{i-\frac{1}{2}j-\frac{1}{2}}}{\Delta y}, \\
	c_{ij} &= h_{i-\frac{1}{2}j-\frac{1}{2}},
\end{flalign}
where $\Delta x$ and $\Delta y$ indicate the grid intervals in the $x$ and $y$ directions, respectively, which are equal to $H (l_{max})$ in this study.

Based on the bilinear representation of the terrain surface, the volumes and areas of the cut-cells are computed. To enable easy computation of the 3D cut-cell parameters, we use the approach of \cite{LockEA2012}, which divides a cell into narrow rectangular columns and approximates the cut-cell volume by the sum of those volumes. The areas of the cut-cells are also computed based on the same approximation. Detailed explanation of this computation is found in \cite{LockEA2012}; see their Appendix B. 

\begin{figure}[t]
  \centering
  \noindent\includegraphics[width=19pc,angle=0]{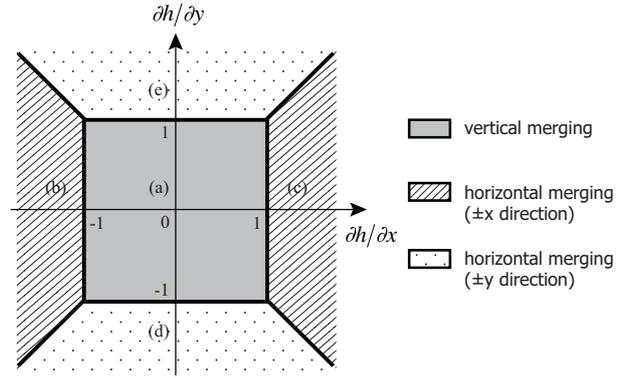}\\
  \caption{The direction of cell-merging depending on the horizontal gradients of the bilinear surface. The shaded region (a) describes the range of the gradients where a small cell is merged with each upper cell. Hatched regions describe the range of the gradients where small cells are merged with each adjacent cell in the (b) $+x$ and (c) $-x$ directions, respectively, and dotted regions describe the same in the (d) $+y$ and (e) $-y$ directions, respectively.}
  \label{Fig5}
\end{figure}

Cell merging is an important feature in the model described in this study. In this method, cut-cells whose center is underground or whose volume is smaller than half the volume of a regular cell are merged with an adjacent cell either vertically or horizontally. The direction of cell-merging is determined by the horizontal gradients of bilinear surfaces that are evaluated at the center of each grid column. In the case of Figure \ref{Fig4}, the gradients of the surface in the $x$ and $y$ directions are evaluated at the grid column ($i, j$) as
\begin{flalign}
	\left(\frac{\partial h}{\partial x}\right)_{ij} &= \frac{(h_{i+\frac{1}{2}j-\frac{1}{2}}+h_{i+\frac{1}{2}j+\frac{1}{2}})-(h_{i-\frac{1}{2}j-\frac{1}{2}}+h_{i-\frac{1}{2}j+\frac{1}{2}})}{2\Delta x}, \label{eq:xgrad} \\
	\left(\frac{\partial h}{\partial y}\right)_{ij} &= \frac{(h_{i-\frac{1}{2}j+\frac{1}{2}}+h_{i+\frac{1}{2}j+\frac{1}{2}})-(h_{i-\frac{1}{2}j-\frac{1}{2}}+h_{i+\frac{1}{2}j-\frac{1}{2}})}{2\Delta y}, \label{eq:ygrad}
\end{flalign}
respectively. When both $|(\partial h/\partial x)_{ij}|$ and $|(\partial h/\partial y)_{ij}|$ are less than or equal to 1, small cut-cells at the column ($i,j$) are merged vertically with each upper cell. Otherwise, they are merged with an adjacent cell in one of four horizontal directions, $+x$, $-x$, $+y$ or $-y$, determined by the values of the gradients (Figure \ref{Fig5}). Note that we switch the direction of cell-merging between the vertical direction and the horizontal directions at the gradients of $\pm 1$.  Following the result of YS10 that the vertical and horizontal merging of cells gave consistent results of flow over a pyramidal mountain, we use the vertical merging at those gradients at our own choice. Similarly, when $|(\partial h/\partial x)_{ij}|$ and $|(\partial h/\partial y)_{ij}|$ are larger than 1 and $|(\partial h/\partial x)_{ij}| = |(\partial h/\partial y)_{ij}|$, we use the horizontal merging in $\pm x$ direction rather than that in $\pm y$ direction. This cell-merging procedure can be described in a form of the pseudocode shown in Figure \ref{Fig6}.

After the cell-merging procedure, the computational volumes of the cells become larger than half of a regular cell, therefore allowing the use of up to half of the full-time step defined by the regular cell. In this study we assume that the topography is well resolved on the grid, and are not concerned with the topography such as an extremely steep v-shaped valley, where a small cut-cell may not have a large adjacent cell to merge with.  

\begin{figure*}
\begin{align*}
&\textit{if} \, \{  \textit{the cell} \:\: (i,j,k) \, \textit{needs to be merged} \, \} \, \textit{then} \nonumber \\
&\quad\qquad\textit{if} \, \{|(\partial h/\partial x)_{ij}| \leqq 1 \:\: \textit{and} \:\:  |(\partial h/\partial y)_{ij}| \leqq 1\}\, \textit{then}\nonumber\\
&\qquad\qquad\qquad \textit{merge with the cell}  \:\: (i,j,k+1) \nonumber\\
&\quad\qquad\textit{elseif}\nonumber\\
&\qquad\qquad\qquad \textit{if} \, \{(\partial h/\partial x)_{ij} > 1 \:\: \textit{and} \:\:  |(\partial h/\partial x)_{ij}| \geqq |(\partial h/\partial y)_{i,j}| \} \,\textit{then}\nonumber\\
&\quad\qquad\qquad\qquad\qquad  \textit{merge with the cell}  \:\: (i-1,j,k)\nonumber\\
&\qquad\qquad\qquad \textit{elseif} \, \{(\partial h/\partial x)_{ij} < -1 \:\: \textit{and} \:\:  |(\partial h/\partial x)_{ij}| \geqq |(\partial h/\partial y)_{i,j}| \} \,\textit{then}\nonumber \\
&\quad\qquad\qquad\qquad\qquad  \textit{merge with the cell} \:\: (i+1,j,k)\nonumber\\
&\qquad\qquad\qquad \textit{elseif} \, \{(\partial h/\partial y)_{ij} > 1 \:\: \textit{and} \:\:  |(\partial h/\partial x)_{ij}| < |(\partial h/\partial y)_{i,j}| \} \, \textit{then}\nonumber\\
&\quad\qquad\qquad\qquad\qquad  \textit{merge with the cell} \:\: (i,j-1,k)\nonumber\\
&\qquad\qquad\qquad \textit{elseif} \, \{(\partial h/\partial y)_{ij} < -1 \:\: \textit{and} \:\:  |(\partial h/\partial x)_{ij}| < |(\partial h/\partial y)_{i,j}| \} \, \textit{then}\nonumber\\
&\quad\qquad\qquad\qquad\qquad  \textit{merge with the cell} \:\: (i,j+1,k)\nonumber\\
&\qquad\qquad\qquad\textit{endif}\nonumber\\
&\quad\qquad\textit{endif}\nonumber\\
&\textit{endif}
\end{align*}
\caption{Pseudocode of the cell-merging algorithm.}
\label{Fig6}
\end{figure*}

\begin{figure*}[t]
  \centering
  \noindent\includegraphics[width=37pc,angle=0]{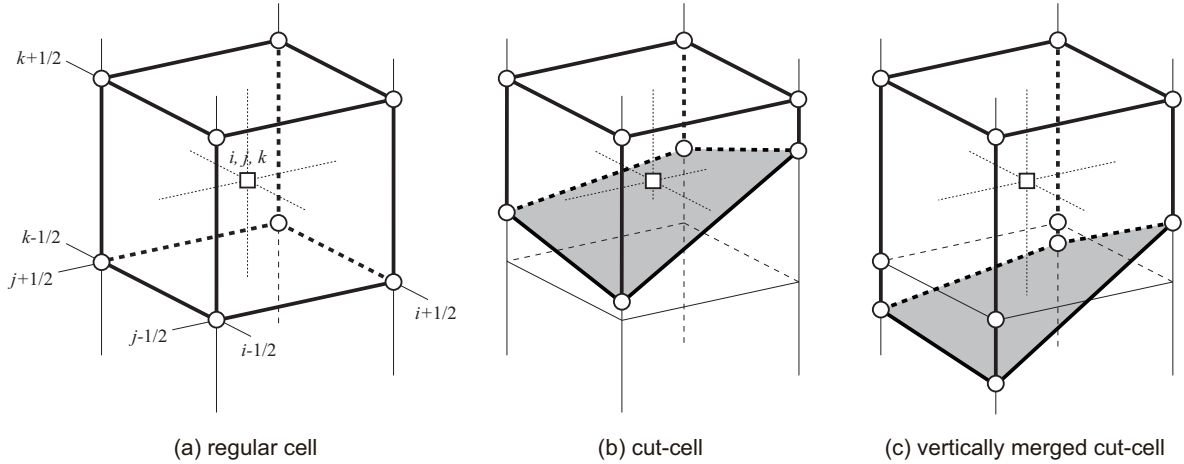}\\
  \caption{Variable arrangement on (a) a regular cell, (b) a non-merged cut cell and (c) a vertically merged cut cell.  Thin solid and dashed lines describe the grid lines, and thick solid and thick-dashed lines describe the boundaries of the cells. Squares and circles represent scalar points and velocity points, respectively. Shaded region represents the topographic surface.}
  \label{Fig7}
\end{figure*}

Finally, the model variables are arranged on the cells. Following the 2D method of YS10, a semi-staggered arrangement of variables is used in this study: scalar variables ($p'$, $\rho'$ and $\theta$) are arranged on the cell centers, and all the velocity components ($u$, $v$ and $w$) are co-located and arranged on the corners of the cells. For an uncut grid cell centered on ($i,j,k$), velocity components are arranged on the eight corners at ($i\pm1/2, j\pm1/2, k\pm1/2$) as shown in Figure \ref{Fig7}(a). Here the square and circles represent the location of the scalar point and velocity points on the cell, respectively. For cut-cells, velocity components are arranged on the corners of the cells above the surface and also on the topographic boundary, as shown for a case with and without cell-merging in Figures \ref{Fig7}(b) and \ref{Fig7}(c), respectively. Variables on horizontally merged cells are also arranged in the same way. Note that no variables are arranged underground. This unique arrangement of variables enables a direct evaluation of the boundary velocity at the terrain surface, as in the 2D method of YS10, thereby simplifying the computation of the velocity near the boundary. The details are discussed in the following subsection as well as the computation of the scalar variables on cut-cells.

\subsection{Spatial discretization}\label{ss:discretization}

To solve flows through the irregularly shaped cut-cells, Eqs (\ref{eq:momentum_u})--(\ref{eq:continuity}) are discretized in space using a finite-volume approach based on the 2D method of YS10. The method invokes Gauss's divergence theorem, which states that the volume integral of the vector divergence over a control volume $V$ enclosed by the surface $S$ is transformed to a surface integral as
\begin{equation}
\int\limits_{V}\nabla\cdot\mathbf{F}\hspace{3pt}dV=\oint\limits_{S}\mathbf{F}\cdot\mathbf{n}\hspace{3pt}dS, \label{eq:gauss}
\end{equation}
where $\mathbf{F}$ is a flux vector and $\mathbf{n}$ is a unit vector pointing along the outward normal of the surface $S$. Assuming that $\nabla\cdot\mathbf{F}$ is constant over a control volume, the vector divergence through a discrete control volume $I$ is expressed as
 \begin{equation}
 (\nabla\cdot\mathbf{F})_{I} = \frac{1}{V_{I}}\oint\limits_{S_{I}}\mathbf{F}\cdot\mathbf{n}\hspace{3pt}dS, \label{eq:gauss_m}
\end{equation}
where $V_{I}$ and $S_{I}$ is the volume and surface of the control volume $I$. Then the surface integral is decomposed into the integrals over all the external sides of the volume. Introducing $F_{I,J}$ as the area mean of the component of $\mathbf{F}$ normal to the side $J$ of the control volume $I$, the surface integral in Eq. (\ref{eq:gauss_m}) becomes 
 \begin{equation}
\oint\limits_{S_{I}}\mathbf{F}\cdot\mathbf{n}\hspace{3pt}dS =   \sum_{J} \int_{S_{I,J}}\mathbf{F}\cdot\mathbf{n}\hspace{3pt}dS  = \sum_{J}F_{I,J}S_{I,J}. \label{eq:gauss_d}
\end{equation}
where $S_{I,J}$ is the area of the side $J$. 

Starting with the simplest case, consider a regular uncut cell shown in Figure \ref{Fig7}(a). Application of Eqs (\ref{eq:gauss_m}) and (\ref{eq:gauss_d}) on the cell ($i,j,k$) gives
\begin{flalign}
(\nabla\cdot\mathbf{F})_{ijk}  &=  \frac{1}{V_{ijk}} \left\{ \delta_{x}(F_{x}S_{x})_{ijk}+\delta_{y}(F_{y}S_{y})_{ijk}+\delta_{z}(F_{z}S_{z})_{ijk} \right\} \nonumber \\
&= \frac{1}{V_{ijk}} \left( F_{x_{i+\frac{1}{2}jk}}S_{x_{i+\frac{1}{2}jk}} - F_{x_{i-\frac{1}{2}jk}}S_{x_{i-\frac{1}{2}jk}} \right. \nonumber \\
&+ \left. F_{y_{ij+\frac{1}{2}k}}S_{y_{ij+\frac{1}{2}k}} - F_{y_{ij-\frac{1}{2}k}}S_{y_{ij-\frac{1}{2}k}} \right. \nonumber \\
&+ \left. F_{z_{ijk+\frac{1}{2}}}S_{z_{ijk+\frac{1}{2}}} - F_{z_{ijk-\frac{1}{2}}}S_{z_{ijk-\frac{1}{2}}} \right),
\label{eq:fvm}
\end{flalign}
where
\begin{eqnarray}
\delta_{x}(\psi)_{ijk} &\equiv& \psi_{i+\frac{1}{2}jk} - \psi_{i-\frac{1}{2}jk} \hspace{1pt},
\end{eqnarray}
and we define $\delta_{y}(\psi)_{ijk}$ and $\delta_{z}(\psi)_{ijk}$ in a similar manner. Here $F_{x}$, $F_{y}$ and $F_{z}$ are the $x$, $y$ and $z$ components of the flux vector, respectively. On a regular cell, the cell volume $V_{ijk}$ and the surface areas normal to the $x$, $y$ and $z$ directions, $S_{x}$, $S_{y}$ and $S_{z}$, respectively, are computed as
\begin{eqnarray}
	V_{ijk} &=& \Delta x \Delta y \Delta z, \label{eq:vol} \\
	S_{x_{i\pm\frac{1}{2}jk}} &=& \Delta y \Delta z, \\
	S_{y_{ij\pm\frac{1}{2}k}} &=& \Delta z \Delta x, \\
	S_{z_{ijk\pm\frac{1}{2}}} &=& \Delta x \Delta y \label{eq:area}, 
\end{eqnarray}
where $\Delta z$ indicates the grid interval in the $z$ direction. Supposing that $\mathbf{F}$ is an advective flux $\mathbf{F} = \phi\mathbf{u}$, where $\phi$ is the scalar quantity, Eq. (\ref{eq:fvm}) becomes
\begin{flalign}
(\nabla\cdot\mathbf{F})_{ijk} &= \frac{1}{V_{ijk}} \left\{ \delta_{x}(\overline{\mathstrut \phi}^{x} \overline{\mathstrut u}^{yz} S_{x})_{ijk} +\delta_{y}(\overline{\mathstrut \phi}^{y} \overline{\mathstrut v}^{zx} S_{y})_{ijk} \right . \nonumber \\
&+ \left. \delta_{z}(\overline{\mathstrut \phi}^{z} \overline{\mathstrut w}^{xy}  S_{z})_{ijk} \right\} \label{fvm_d}, 
\end{flalign}
where
\begin{eqnarray}
\overline{\psi}_{ijk}^{x} &\equiv& \left(\psi_{i-\frac{1}{2}jk}+\psi_{i+\frac{1}{2}jk}\right)/2 \hspace{1pt},
\end{eqnarray}
and we define $\overline{\psi}_{ijk}^{y}$ and $\overline{\psi}_{ijk}^{z}$ in a similar manner. When the volume and areas of a regular cell (\ref{eq:vol})--(\ref{eq:area}) are assigned, Eq. (\ref{fvm_d}) reduces to a centered finite-difference expression for $\nabla\cdot\mathbf{F}$. 
 
Using the notations introduced above, the spatially discretized forms of Eqs (\ref{eq:momentum_u})--(\ref{eq:continuity}) for an uncut cell are given as follows: 
\begin{flalign}
\left(\frac{\partial\rho u}{\partial t}\right)_{i'j'k'} &= - \frac{1}{V_{i'j'k'}} \left\{ \delta_{x}(\overline{\mathstrut \rho}^{yz}\overline{\mathstrut u}^{x}\overline{\mathstrut u}^{x}S_{x})_{i'j'k'} \right. \nonumber \\
&+ \left. \delta_{y}(\overline{\mathstrut \rho}^{zx}\overline{\mathstrut u}^{y}\overline{\mathstrut v}^{y}S_{y})_{i'j'k'} 
+ \delta_{z}(\overline{\mathstrut \rho}^{xy}\overline{\mathstrut u}^{z}\overline{\mathstrut w}^{z}S_{z})_{i'j'k'} \right\} \nonumber \\
&- \frac{\delta_{x}(\overline{\mathstrut p^{\prime}}^{yz})_{i'j'k'}}{\Delta x} + {D_{u}^{*}}_{i'j'k'},   \label{eq:momentum_u_d} \\
\left(\frac{\partial\rho v}{\partial t}\right)_{i'j'k'} &= - \frac{1}{V_{i'j'k'}} \left\{ \delta_{x}(\overline{\mathstrut \rho}^{yz}\overline{\mathstrut v}^{x}\overline{\mathstrut u}^{x}S_{x})_{i'j'k'} \right. \nonumber\\
&+ \left. \delta_{y}(\overline{\mathstrut \rho}^{zx}\overline{\mathstrut v}^{y}\overline{\mathstrut v}^{y}S_{y})_{i'j'k'} 
+ \delta_{z}(\overline{\mathstrut \rho}^{xy}\overline{\mathstrut v}^{z}\overline{\mathstrut w}^{z}S_{z})_{i'j'k'} \right\}  \nonumber \\
&- \frac{\delta_{y}(\overline{p^{\prime}}^{zx})_{i'j'k'}}{\Delta y} + {D_{v}^{*}}_{i'j'k'},   \label{eq:momentum_v_d} \\
\left(\frac{\partial\rho w}{\partial t}\right)_{i'j'k'}  &= - \frac{1}{V_{i'j'k'}}\left\{\delta_{x}(\overline{\rho}^{yz}\overline{w}^{x}\overline{u}^{x}S_{x})_{i'j'k'} \right. \nonumber \\
 &+ \left. \delta_{y}(\overline{\rho}^{zx}\overline{w}^{y}\overline{v}^{y}S_{y})_{i'j'k'} 
+ \delta_{z}(\overline{\rho}^{xy}\overline{w}^{z}\overline{w}^{z}S_{z})_{i'j'k'}\right\} \nonumber \\
&- \frac{\delta_{z}(\overline{p^{\prime}}^{xy})_{i'j'k'}}{\Delta z}  - \overline{\rho^{\prime}}^{xyz}_{i'j'k'}g +{D_{w}^{*}}_{i'j'k'},  \label{eq:momentum_w_d} \\ 
\left(\frac{\partial p^{\prime}}{\partial t}\right)_{ijk}  &=  - \frac{c_{p}R}{c_{v}}\left(\frac{p_{ijk}}{p_{0}}\right)^{\frac{R}{c_{p}}} \frac{1}{V_{ijk}}\left\{\delta_{x}(\overline{\mathstrut \rho\theta}^{x}\overline{\mathstrut u}^{yz}S_{x})_{ijk} \right . \nonumber \\
&+ \left. \delta_{y}(\overline{\mathstrut \rho\theta}^{y}\overline{\mathstrut v}^{zx}S_{y})_{ijk}
+ \delta_{z}(\overline{\mathstrut \rho\theta}^{z}\overline{\mathstrut w}^{xy}S_{z})_{ijk} \right. \nonumber \\
&- \left. {D_{\theta}^{*}}_{ijk}\right\}, \\ 
\left(\frac{\partial\rho^{\prime}}{\partial t}\right)_{ijk}  &= - \frac{1}{V_{ijk}}\left\{\delta_{x}(\overline{\rho}^{x}\overline{u}^{yz}S_{x})_{ijk} + \delta_{y}(\overline{\rho}^{y}\overline{v}^{zx}S_{y})_{ijk} \right. \nonumber \\
&+ \left. \delta_{z}(\overline{\rho}^{z}\overline{w}^{xy}S_{z})_{ijk}\right\}, 
\end{flalign}
where $V_{i'j'k'}$ is the volume of a velocity cell centered on ($i',j',k'$) = ($i-1/2, j-1/2, k-1/2$), which is equal to $V_{ijk}$ when the velocity cell is uncut. The terms $D_{u}^{*}$, $D_{v}^{*}$, $D_{w}^{*}$, and $D_{\theta}^{*}$ are discretized forms of the diffusion terms. In this study the second-order central difference scheme is used for pressure gradient terms and diffusion terms. For cells without necessary neighboring fluid points for the calculation of a fourth-order artificial diffusion term, a second-order term is introduced instead, where boundary conditions are used at the terrain surface. 

\begin{figure}[t]
  \centering
  \noindent\includegraphics[width=19pc,angle=0]{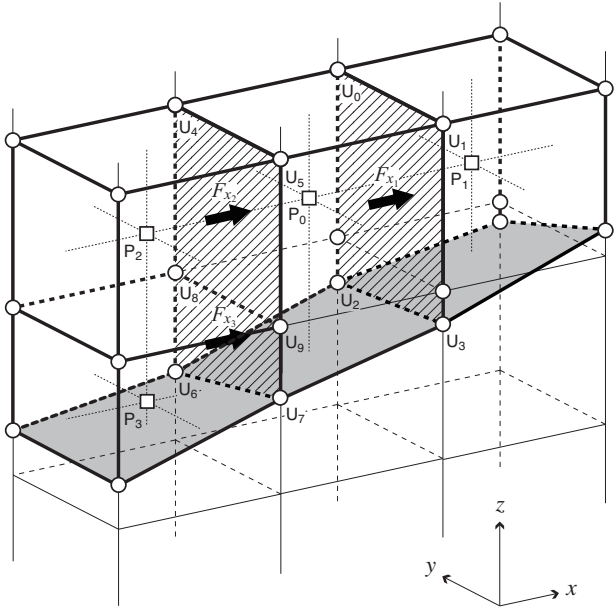}\\
  \caption{Flux calculation on a scalar cell cut by the topography. Arrows indicate fluxes through the faces normal to the $x$ direction of the scalar cell centered at $\mathrm{P_{0}}$. Hatched regions indicate the cell faces normal to the $x$ direction. Thin solid and dashed lines describe the grid lines, and thick solid and thick-dashed lines describe the boundaries of the cells. Squares and circles represent scalar points and velocity points, respectively. Shaded region represents the topographic surface.  }
  \label{Fig8}
\end{figure}

We will now describe how the finite-volume approach is implemented in a situation where some of the cells are cut by the topography. As a semi-staggered arrangement of the scalar and velocity variables is used, the scalar and velocity cells are at different locations and are treated separately. First, consider the discretization of conservation equations on the scalar cell centered at point $\mathrm{P_{0}}$ (Figure \ref{Fig8}). The cell is vertically merged, and enclosed by eight velocity points at $\mathrm{U}_{0}$ to $\mathrm{U}_{7}$. Though we focus here on the case of a vertically merged cell, the spatial discretization on horizontally merged cells is handled in the same way. 

As shown in the case of an uncut cell, finite-volume discretization requires the estimation of surface integrals over each face of the cell. For example, the surface integral of the flux over the right face of the cell $\mathrm{P_{0}}$ is evaluated as 
\begin{equation}
\int_{\mathrm{U}_{0123}}\mathbf{F}\cdot\mathbf{n}\hspace{3pt}dS = F_{x_{1}}S_{x_{1}},
\end{equation}
where $\mathrm{U}_{0123}$ denotes the face enclosed by the points $\mathrm{U}_{0}$, $\mathrm{U}_{1}$, $\mathrm{U}_{2}$ and $\mathrm{U}_{3}$. In the right-hand side, $S_{x_{1}}$ represents the area of the face $\mathrm{U_{0123}}$, and $F_{x_{1}}$ is the area mean of the flux over the face (Figure \ref{Fig8}). In the case of the left face of the cell $\mathrm{P_{0}}$, $\mathrm{U}_{4567}$, the flux through the face is composed of two fluxes: the flux through the boundary with the cell $\mathrm{P_{2}}$ and that through the boundary with the cell $\mathrm{P_{3}}$. Thus the surface integral of the flux over the face $\mathrm{U}_{4567}$ can be decomposed as 
\begin{flalign}
\int_{\mathrm{U}_{4567}}\mathbf{F}\cdot\mathbf{n}\hspace{3pt}dS &= \int_{\mathrm{U}_{4589}}\mathbf{F}\cdot\mathbf{n}\hspace{3pt}dS + \int_{\mathrm{U}_{6789}}\mathbf{F}\cdot\mathbf{n}\hspace{3pt}dS \\
&= F_{x_{2}}S_{x_{2}}+F_{x_{3}}S_{x_{3}}, 
\end{flalign}
where $S_{x_{2}}$ and $S_{x_{3}}$ are the areas of the boundary faces $\mathrm{U_{4589}}$ and $\mathrm{U_{6789}}$, respectively, and $F_{x_{2}}$ and $F_{x_{3}}$ are the area means of flux over each of those faces, respectively.

The evaluation of $F_{x_{1}}$, $F_{x_{2}}$ and $F_{x_{3}}$ requires the area mean of the normal velocity and the advected scalar quantity over each face. In this study, the normal velocity of these fluxes is estimated by a linear interpolation among the velocity values at the corners of each face. For example, the normal velocity of $F_{x_{1}}$ is computed as 
\begin{eqnarray}
u_{x_{1}} &=& (u_{_{\mathrm{U}_{0}}}+u_{_{\mathrm{U}_{1}}}+u_{_{\mathrm{U}_{2}}}+u_{_{\mathrm{U}_{3}}})/4. \label{eq:ux1} 
\end{eqnarray}
The normal velocity of $F_{x_{2}}$ and $F_{x_{3}}$ on the merged face $\mathrm{\mathrm{U}}_{4567}$ are given by
\begin{eqnarray}
u_{x_{2}} &=& (u_{_{\mathrm{U}_{4}}}+u_{_{\mathrm{U}_{5}}}+u_{_{\mathrm{U}_{8}}}+u_{_{\mathrm{U}_{9}}})/4, \label{eq:ux2}  \\
u_{x_{3}} &=& (u_{_{\mathrm{U}_{6}}}+u_{_{\mathrm{U}_{7}}}+u_{_{\mathrm{U}_{8}}}+u_{_{\mathrm{U}_{9}}})/4, \label{eq:ux3} 
\end{eqnarray}
respectively. Note that the calculation of the velocity values in the right-hand sides of the Eqs (\ref{eq:ux1})--(\ref{eq:ux3}) is described later in this subsection. 

The calculation of the area mean of the scalar quantity over the faces can be more complicated. In the case of Figure \ref{Fig8}, the scalar quantity in $F_{x_{2}}$ is obtained by a simple linear interpolation between neighboring cell centers as
\begin{equation}
\phi_{x_{2}} = (\phi_{_{\mathrm{P}_{0}}}+\phi_{_{\mathrm{P}_{2}}})/2.
\end{equation}
The evaluation of the scalar quantity in $F_{x_{3}}$, on the other hand, is not straightforward because one of the neighboring cell centers is underground and absorbed to the cell $\mathrm{P_{0}}$ in the cell-merging procedure. One way to evaluate it with a high-order accuracy is to use a multi-dimensional interpolation or extrapolation scheme around the boundary. A disadvantage of the scheme is that it generally entails a considerable complexity, especially in 3D, to locate the appropriate grid points to form a multi-dimensional function for various shapes of cut-cells. For computational simplicity, here we use a first-order calculation as used in the 2D study of YS12. 
In this method, the area mean of the scalar quantity over a face is approximated by the simple average of cell center values of the cells exchanging fluxes through the face. In the case of Figure \ref{Fig8}, the scalar quantity in $F_{x_{1}}$ and $F_{x_{3}}$ is therefore evaluated as
\begin{flalign}
\phi_{x_{1}} &= (\phi_{_{\mathrm{P}_{0}}}+\phi_{_{\mathrm{P}_{1}}})/2, \\
\phi_{x_{3}} &= (\phi_{_{\mathrm{P}_{0}}}+\phi_{_{\mathrm{P}_{3}}})/2,
\end{flalign}
respectively. Note that this approximation does not violate mass-conservation. As a result, the surface integrals of the flux over the left and right faces of the cell $\mathrm{P}_{0}$ are computed as 
\begin{flalign}
\int_{\mathrm{U}_{4567}}\mathbf{F}\cdot\mathbf{n}\hspace{3pt}dS &= \phi_{x_{2}}u_{x_{2}}S_{x_{2}} + \phi_{x_{3}}u_{x_{3}}S_{x_{3}}, \\
\int_{\mathrm{U}_{0123}}\mathbf{F}\cdot\mathbf{n}\hspace{3pt}dS &= \phi_{x_{1}}u_{x_{1}}S_{x_{1}},
\end{flalign}
respectively. The surface integrals of the flux over the other faces of the cell, $\mathrm{U}_{1357}$, $\mathrm{U}_{0246}$ and $\mathrm{U}_{0145}$, are also computed in the same way, with zero normal flux assumed at the topographic boundary. Boundary conditions that involve nonzero normal fluxes, for example heat conduction across the surface, can be written as effective volume-mean source terms \citep[][]{AdcroftEA1997}. 

Next we consider the discretization of momentum equations on the velocity cells. A key problem here is to accurately evaluate the pressure gradient at each velocity cell. Because some pressure points are underground in a cut-cell model, a velocity cell cut by the topography may not have all the necessary pressure points for the calculation of the pressure gradient. To allow the same pressure gradient calculation on cut-cells as is used for a regular cell, \cite{Walko+Avissar2008} predicted approximate solutions of underground pressure by assuming that the cell volumes and cell areas of cut-cells are uniformly occupied by small solids. The assumption makes the shape of topography indistinct and thus could affect the flow dynamics near the topography. \cite{YeEA1999}, on the other hand, proposed an approach to express the pressure field near the surface in terms of a polynomial interpolating function, and evaluate the gradients on cut-cells based on the interpolating function. They used pressure values at available neighboring cell-centers and also at the surface to construct an interpolating function for each cut-cell. This method allows the systematic evaluation of the pressure gradient at cut-cells of various shapes. However the calculation of surface pressure causes difficulties in atmospheric simulations because the zero-gradient boundary condition on the pressure field does not hold true for density-stratified flow. In addition, since geophysical flow is driven by slight perturbations of the pressure from the hydrostatically balanced state, extrapolation of the pressure perturbation to the surface could also affect the local flow dynamics as well as use of the underground pressure. 

The commonly-used staggered arrangement of variables in atmospheric models introduces another complexity into the discretization of momentum equations. On a non-staggered grid, not only are the velocity components and scalar variables co-located, but the position and geometry of the associated cells are also identical. With a staggered grid, the scalar and velocity cells are at different locations and will generally have a different shape when they are cut by the topography. In general, a cut cell method for a staggered grid must deal with this extra complexity \citep[][]{KirkpatrickEA2003}. With a cell-merging approach, in particular, another formulation of the discretized equations for merged velocity cells is normally required. 

\begin{figure}[t]
  \centering
  \noindent\includegraphics[width=19pc,angle=0]{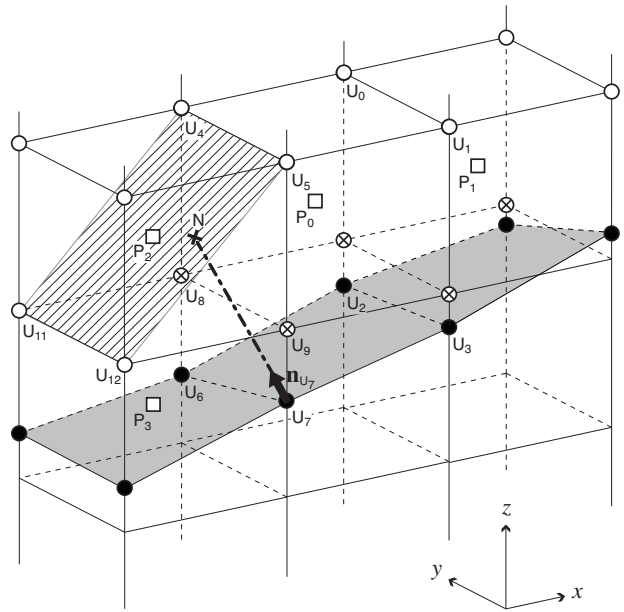}\\
  \caption{Velocity calculation on cut-cells. Squares represent scalar points. Open circles represent velocity points on which model solutions are predicted. Filled circles and crossed circles represent diagnosed velocity points on the terrain surface and on the merged faces of the scalar cells, respectively. The dashed-dotted line and the arrow describe the normal line and the unit normal direction to the surface at the boundary point $\mathrm{U}_{7}$, respectively. The hatched region indicates the nearest plane for the boundary point $\mathrm{U}_{7}$ defined by four open-circled points of the scalar cell $\mathrm{P}_{2}$.}
  \label{Fig9}
\end{figure}

The current method provides a way around these complications by using the non-conventional semi-staggered variable arrangement described in section \ref{ss:cut-cells}. For example, consider the velocity points illustrated in Figure \ref{Fig9}, which are arranged on the same scalar cells shown in Figure \ref{Fig8}. First, we solve the momentum equations on the velocity cells which retain the regular rectangular shape. Supposing that the pressure points neighboring $\mathrm{P}_{0}$, $\mathrm{P}_{1}$, $\mathrm{P}_{2}$ and $\mathrm{P}_{3}$ are all available, the velocity cells with centers indicated by open circles have all the eight pressure values available at the corners of the cells. Therefore the velocity values on those cells are obtained using Eqs (\ref{eq:momentum_u_d})--(\ref{eq:momentum_w_d}). Once the velocity at the open-circled points is predicted, the remaining velocity points are either on the topographic surface or on the merged faces of the scalar cells, indicated by filled circles and crossed circles, respectively (Figure \ref{Fig9}). The velocity on these points are evaluated diagnostically. 

The velocity on the topographic surface is calculated by applying a velocity boundary condition on the surface. When the non-slip boundary condition is imposed, all the boundary velocities on the surface are set to zero. With the free-slip boundary condition, the boundary velocity is calculated so that the component of the velocity that is tangential to the surface is preserved near the boundary. Given that the unit normal direction to the surface at a boundary point $\mathrm{U}_{\mathrm{B}}$ is $\mathbf{n}_{_{\mathrm{U}_{\mathrm{B}}}}$, the boundary velocity $\mathbf{u}_{_{\mathrm{U}_{\mathrm{B}}}}$ is therefore calculated as
\begin{eqnarray}
 	\mathbf{u}_{_{\mathrm{U}_{\mathrm{B}}}} = \tilde{\mathbf{u}}_{_{\mathrm{U}_{\mathrm{B}}}} -(\tilde{\mathbf{u}}_{_{\mathrm{U}_{\mathrm{B}}}}\cdot\mathbf{n}_{_{\mathrm{U}_{\mathrm{B}}}})\mathbf{n}_{_{\mathrm{U}_{\mathrm{B}}}} \hspace{3pt}, \label{eq:33}
\end{eqnarray}
where $\tilde{\mathbf{u}}_{_{\mathrm{U}_{\mathrm{B}}}}$ indicates a mean velocity near the boundary point $\mathrm{U}_{\mathrm{B}}$ that is computed from the predicted velocities on the regular velocity cells.

For example, consider the computation of the velocity at the boundary point $\mathrm{U}_{7}$ using the free-slip condition. Supposing that the scalar cell $\mathrm{P}_{0}$ is at the grid column ($i, j$), the unit normal direction to the surface at the boundary point $\mathrm{U}_{7}$,
\begin{eqnarray}
\mathbf{n}_{_{\mathrm{U}_{7}}} = (\mathrm{n}_{x}, \mathrm{n}_{y}, \mathrm{n}_{z}),
\end{eqnarray}
is computed as,
\begin{eqnarray}
\mathrm{n}_{x} = \frac{-(\partial h/\partial x)_{i-\frac{1}{2}j-\frac{1}{2}}}{\sqrt{(\partial h/\partial x)_{i-\frac{1}{2}j-\frac{1}{2}}^{2}+(\partial h/\partial y)_{i-\frac{1}{2}j-\frac{1}{2}}^{2}+1}}, \\
\mathrm{n}_{y} = \frac{-(\partial h/\partial y)_{i-\frac{1}{2}j-\frac{1}{2}}}{\sqrt{(\partial h/\partial x)_{i-\frac{1}{2}j-\frac{1}{2}}^{2}+(\partial h/\partial y)_{i-\frac{1}{2}j-\frac{1}{2}}^{2}+1}}, \\
\mathrm{n}_{z} = \frac{1}{\sqrt{(\partial h/\partial x)_{i-\frac{1}{2}j-\frac{1}{2}}^{2}+(\partial h/\partial y)_{i-\frac{1}{2}j-\frac{1}{2}}^{2}+1}},
\end{eqnarray}
where the gradients of the surface at the boundary point $\mathrm{U}_{7}$ are calculated at their horizontal grid positions ($i-1, j-1$). Following Eqs (\ref{eq:xgrad}) and (\ref{eq:ygrad}), the gradients are obtained as
\begin{flalign}
	\left(\frac{\partial h}{\partial x}\right)_{i-\frac{1}{2}j-\frac{1}{2}} = \frac{(h_{ij-1}+h_{ij})-(h_{i-1j-1}+h_{i-1j})}{2\Delta x}, \\
	\left(\frac{\partial h}{\partial y}\right)_{i-\frac{1}{2}j-\frac{1}{2}} = \frac{(h_{i-1j}+h_{ij})-(h_{i-1j-1}+h_{ij-1})}{2\Delta y}.
\end{flalign}
Next, to estimate a mean velocity near the boundary point, we define the nearest plane to the point by four open-circled points within a cell. The mean velocity near the point is then evaluated at the intersection of the nearest plane and the normal line to the surface. When the normal line at the boundary point $\mathrm{U}_{7}$ is as described in Figure \ref{Fig9}, the velocities at the four points $\mathrm{U}_{4}$,  $\mathrm{U}_{5}$, $\mathrm{U}_{11}$ and $\mathrm{U}_{12}$ are used to evaluate the mean velocity at the intersection point $\mathrm{N}$ as
\begin{flalign}
 \mathbf{u}_{_{\mathrm{N}}} &=  \alpha\beta \,\mathbf{u}_{_{\mathrm{U}_{4}}}+\alpha(1-\beta)\mathbf{u}_{_{\mathrm{U}_{5}}} + (1-\alpha)\beta\,\mathbf{u}_{_{\mathrm{U}_{11}}}  \nonumber \\
&+ (1-\alpha)(1-\beta)\mathbf{u}_{_{\mathrm{U}_{12}}}, 
\end{flalign} 
where $\alpha$ and $\beta$ are the linear interpolation factors defined as
\begin{eqnarray}
	\alpha = \frac{x_{_{\mathrm{N}}}-x_{_{\mathrm{U_{12}}}}}{\Delta x}, \hspace{3pt}
	\beta = \frac{y_{_{\mathrm{N}}}-y_{_{\mathrm{U_{12}}}}}{\Delta y}. 
\end{eqnarray} 
As a result, the velocity at the boundary point $\mathrm{U}_{7}$ is computed as
\begin{eqnarray}
	\mathbf{u}_{_{\mathrm{U}_{7}}} = \mathbf{u}_{_{\mathrm{N}}} -(\mathbf{u}_{_{\mathrm{N}}}\cdot\mathbf{n}_{_{\mathrm{U}_{7}}})\mathbf{n}_{_{\mathrm{U}_{7}}}. 
\end{eqnarray}
The velocities at the other boundary points, such as at $\mathrm{U}_{2}$, $\mathrm{U}_{3}$ and $\mathrm{U}_{6}$, are also calculated in the same way. 

Finally, the velocities at crossed-circled points on the merged faces are calculated assuming a linear distribution of the velocity over each merged face. For example, the velocity at $\mathrm{U}_{8}$ and $\mathrm{U}_{9}$ on the merged face $\mathrm{U}_{4567}$ of the scalar cell $\mathrm{P}_{0}$ is calculated by a linear interpolation between the velocities at the corners of the merged face as 
\begin{flalign}
	\mathbf{u}_{_{\mathrm{U}_{8}}} &= \mathbf{u}_{_{\mathrm{U}_{4}}}\gamma+\mathbf{u}_{_{\mathrm{U}_{6}}}(1-\gamma) \hspace{1pt}, \\
	\mathbf{u}_{_{\mathrm{U}_{9}}} &= \mathbf{u}_{_{\mathrm{U}_{5}}}\zeta+\mathbf{u}_{_{\mathrm{U}_{7}}}(1-\zeta) \hspace{1pt}, 
\end{flalign}
respectively, where the interpolation factors $\gamma$ and $\zeta$ are defined as
\begin{eqnarray}
	\gamma = \frac{z_{_{\mathrm{U}_{8}}}-z_{_{\mathrm{U_{6}}}}}{z_{_{\mathrm{U}_{4}}}-z_{_{\mathrm{U_{6}}}}}, \hspace{5pt}
	\zeta = \frac{z_{_{\mathrm{U}_{9}}}-z_{_{\mathrm{U_{7}}}}}{z_{_{\mathrm{U}_{5}}}-z_{_{\mathrm{U_{7}}}}}. 
\end{eqnarray}

All the velocity values are calculated without estimating the surface or the underground pressure with the following three-step calculation of the velocity: 1) solve the momentum equations on regular velocity cells, 2) calculate the velocity on the terrain surface using a boundary condition, and 3) calculate the velocity on merged faces by an interpolation of the velocity obtained in 1) and 2). In addition, there is no need to merge velocity cells, and thus the computational cost of merging cells is kept as low as that of non-staggered models. Because the momentum equations are only solved on the fully rectangular cells that do not enclose the entire fluid domain, momentum is not conserved in our model as it is in many other models using a vertically staggered grid (e.g., Saito et al. 2001; Satoh 2002; Klemp et al. 2007). Our model also does not conserve kinetic energy, although mass-conservation is guaranteed. 

The cut-cell calculation described in this subsection is only performed in the finest cubes. The ghost cells that are added beyond the boundary of each cube allow merging of cells between adjacent finest cubes. To prevent cell merging between cells that belong to different sizes of cubes, the size of the cubes near the terrain surface is readjusted in the grid-generation process, if necessary, so that both merged and non-merged cut-cells are inside of the finest region. In the coarser cubes, the standard finite-difference code is used so that the computational overhead induced by the use of cut-cell code only occurs in the area near the terrain surface.  While on cut-cells the scheme is locally first-order, a global second-order accuracy is achieved, as demonstrated numerically in the appendix. 

\subsection{Parallelization}\label{ss:parallel}

The parallelization of the model is straightforward. Load balancing is achieved by distributing an equal number of cubes at each refinement level for each processor. To keep the load balance in the presence of cut-cells and underground cells at the finest level, we compute solutions for all cells in the finest cubes using the same cut-cell code, and then set the values at cells that are located below the terrain surface to zero. The overhead required to calculate values of underground cells is sufficiently low because cells in cubes that are completely bellow the terrain surface are excluded from the computation (Figure \ref{Fig1}c). For efficient parallelization, a sufficient number of cubes would be required to avoid inequality in the distributed numbers of cubes among processors.

The model code is parallelized using OpenMP. The scalability of the model is demonstrated in section \ref{s:refine} through mountain-wave simulations both on a uniform grid and a locally refined grid. The data structure of the block-structured mesh used in this study can also be suitable for parallel computing by Message Passing Interface (MPI) in a distributed memory system, as demonstrated using BCM by \cite{TakahashiEA2008} and \cite{SakaiEA2013}.

\section{Results}\label{s:result}
This section presents the model results from test simulations of flow passing over a 3D hill. Every test involves an isolated hill located at the center of the lower boundary. The model is integrated for 1 hour for each test. The following atmospheric and boundary conditions are imposed on all the simulations. A constant horizontal velocity and Brunt-V\"ais\"al\"a frequency, $U$ = 10 m s$^{-1}$ and $N$ = 0.01 s$^{-1}$, respectively, are initially imposed on the entire domain. A sea level potential temperature is specified to be $\Theta = 300$ K. The lower and lateral boundary conditions are free-slip and cyclic, respectively. To prevent cyclic lateral boundaries from contaminating the simulated results, a large domain length of 64 km is used in the $x$ direction. In the $y$ direction, the domain length is set to 32 km. The height of the domain is 16 km, and a sponge layer \citep{Klemp+Lilly1978} is placed higher than 10 km to avoid reflecting the gravity wave at the rigid top boundary. 

\subsection{Cut-cell model on a uniform mesh}\label{s:uniform}

First, we test the model on a uniform grid of equal spacing in order to validate the cut-cell representation of topography without using any grid refinement. The performance of the model on a locally refined grid is discussed in section \ref{s:refine}.

\subsubsection{Flow over a 3D bell-shaped hill} \label{ss:bell}

This test examines the accuracy of the cut-cell representation of topography by simulating a flow over a 3D bell-shaped hill and comparing the result to that predicted by theory. The surface height of the hill is described as, 
\begin{eqnarray}
h(x,y) = \frac{h_{\mathrm{m}}}{(1+x^{2}/a^{2}+y^{2}/a^{2})^{3/2}},
\end{eqnarray}
where the maximum height $h_{\mathrm{m}}$ = 400 m and the half-width $a$ = 1 km. Following \cite{LockEA2012}, the general shape of the hill is chosen to agree with that analyzed in \cite{Smith1980}. Here the same grid spacing of 200 m is used in x, y and z directions over the entire domain.  

The maximum slope angle of the hill resolved by the model grid is 18.1$^{\circ}$. The horizontal gradients of the terrain surface, calculated by Eqs (\ref{eq:xgrad}) and (\ref{eq:ygrad}), are less than 1 at all grid columns, so that only the vertical merging of cells are used in this test. 

\begin{figure}[t]
  \centering
  \noindent\includegraphics[width=19pc,angle=0]{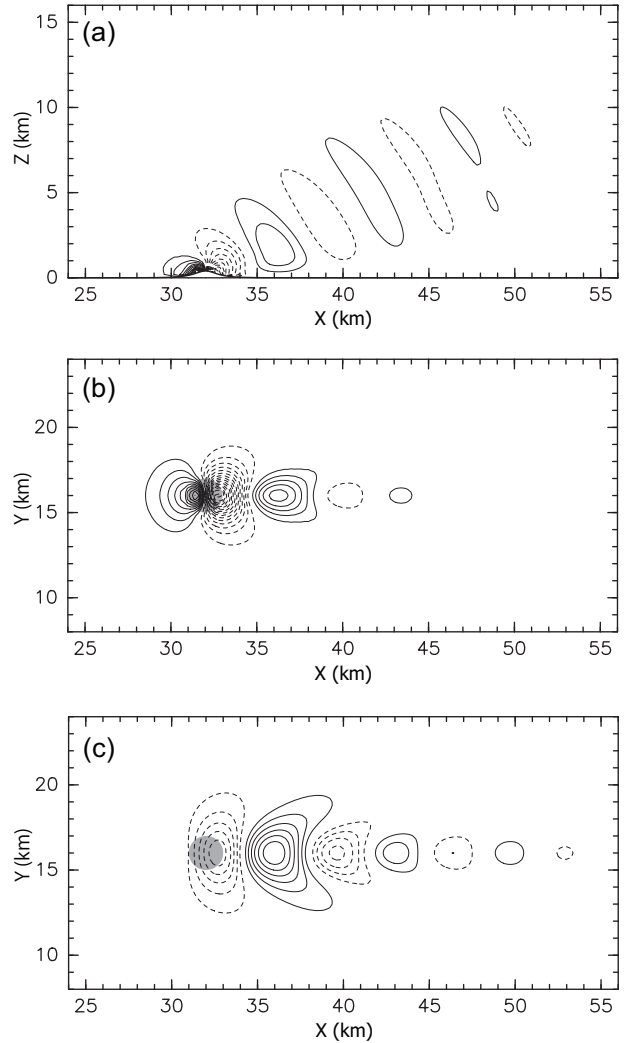}\\
  \caption{The vertical velocity over the 3D bell-shaped hill after 1 hour of integration: (a) a vertical $x$-$z$ slice through the center of the hill, (b) a horizontal $x$-$y$ slice at height $z$ = 800 m and (c) a horizontal $x$-$y$ slice at height $z$ = 2000 m. Contour intervals are 0.25 m s$^{-1}$ in (a) and 0.1 m s$^{-1}$ in (b) and (c). Solid and dashed lines indicate positive and negative values, respectively. Shaded regions in (b) and (c) indicate the region within the half-width from the top of the hill.}
  \label{Fig10}
\end{figure}

Figure \ref{Fig10} shows the simulated vertical velocity fields over the 3D bell-shaped hill. The results successfully demonstrate stable and smooth solutions over many integrations. A vertical $x$-$z$ cross section through the center of the hill (Figure \ref{Fig10}(a)) shows the generation of vertically propagating waves from the top of the hill. The phase lines tilt upstream with increasing height above the terrain surface, just as they do in the 2D flow over a bell-shaped hill with the same height and half-width, such as examined in \cite{Gallus+Klemp2000}. Unlike the 2D case, the amplitude of the flow decays rapidly with height because of the rapid reduction of the energy density due to the three-dimensionality. Figures \ref{Fig10}(b) and (c) show the horizontal $x$-$y$ slices of the vertical velocity field. Above the surface at $z$ = 800 m (Figure \ref{Fig10}(b)), the fluid rises on the upstream side of the hill and descends over the leeward, showing a damped oscillation towards downstream. The wavelength of the oscillation is about 7 km, which is comparable to that calculated by the linear theory in \cite{Smith1980}, $2\pi U/N \approx 6.3$ km, and also agrees well with the numerical result by \cite{LockEA2012}. At the higher altitude of $z$ = 2000 m (Figure \ref{Fig10}(c)), the flow patterns on the leeward of the mountain splits and widens to form U-shaped regions in good agreement with the theory.

\begin{figure}[t]
  \centering
  \noindent\includegraphics[width=19pc,angle=0]{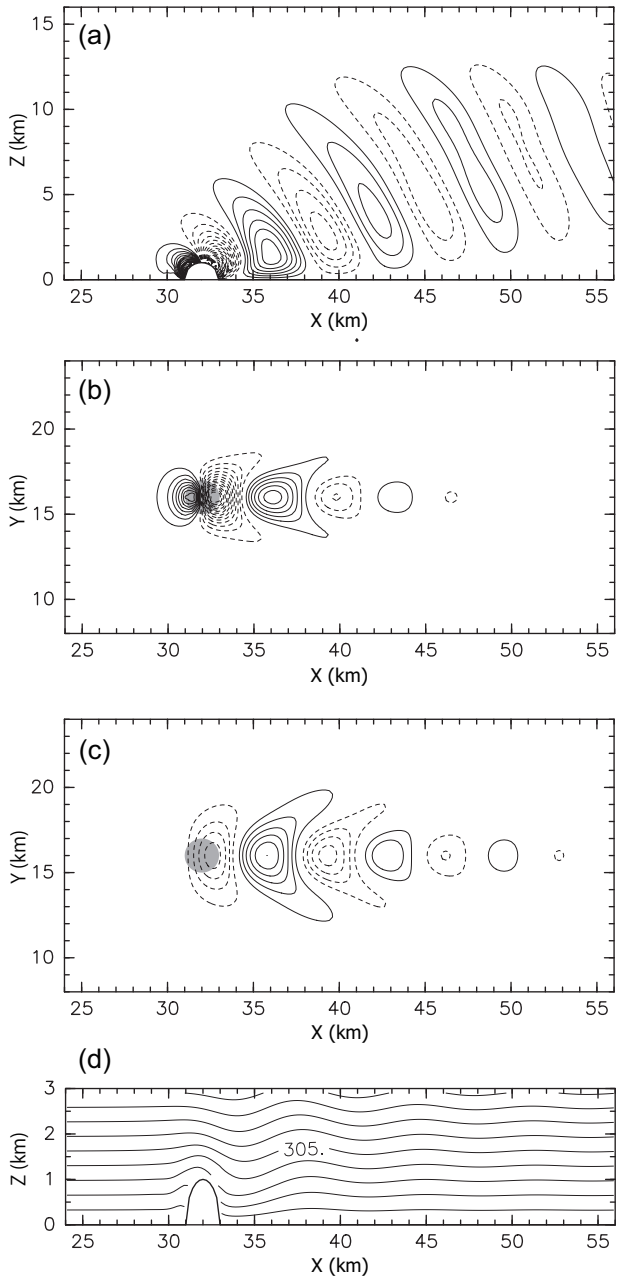}\\
  \caption{Results of flow over the 3D hemisphere hill after 1 hour of integration. (a) Vertical $x$-$z$ slice of the vertical velocity field through the center of the hill. (b)--(c) Horizontal $x$-$y$ slices of the vertical velocity field at height $z$ = 1200 m and $z$ = 2400 m, respectively. (d) Vertical $x$-$z$ slice of the potential temperature field through the center of the hill. Contour intervals are 0.25 m s$^{-1}$ in (a), (b) and (c), and 1 K in (d). Solid and dashed lines indicate positive and negative values, respectively. Shaded regions in (b) and (c) indicate the position of the hill.}
  \label{Fig11}
\end{figure}

\begin{figure*}[t]
  \centering
  \noindent\includegraphics[width=39pc,angle=0]{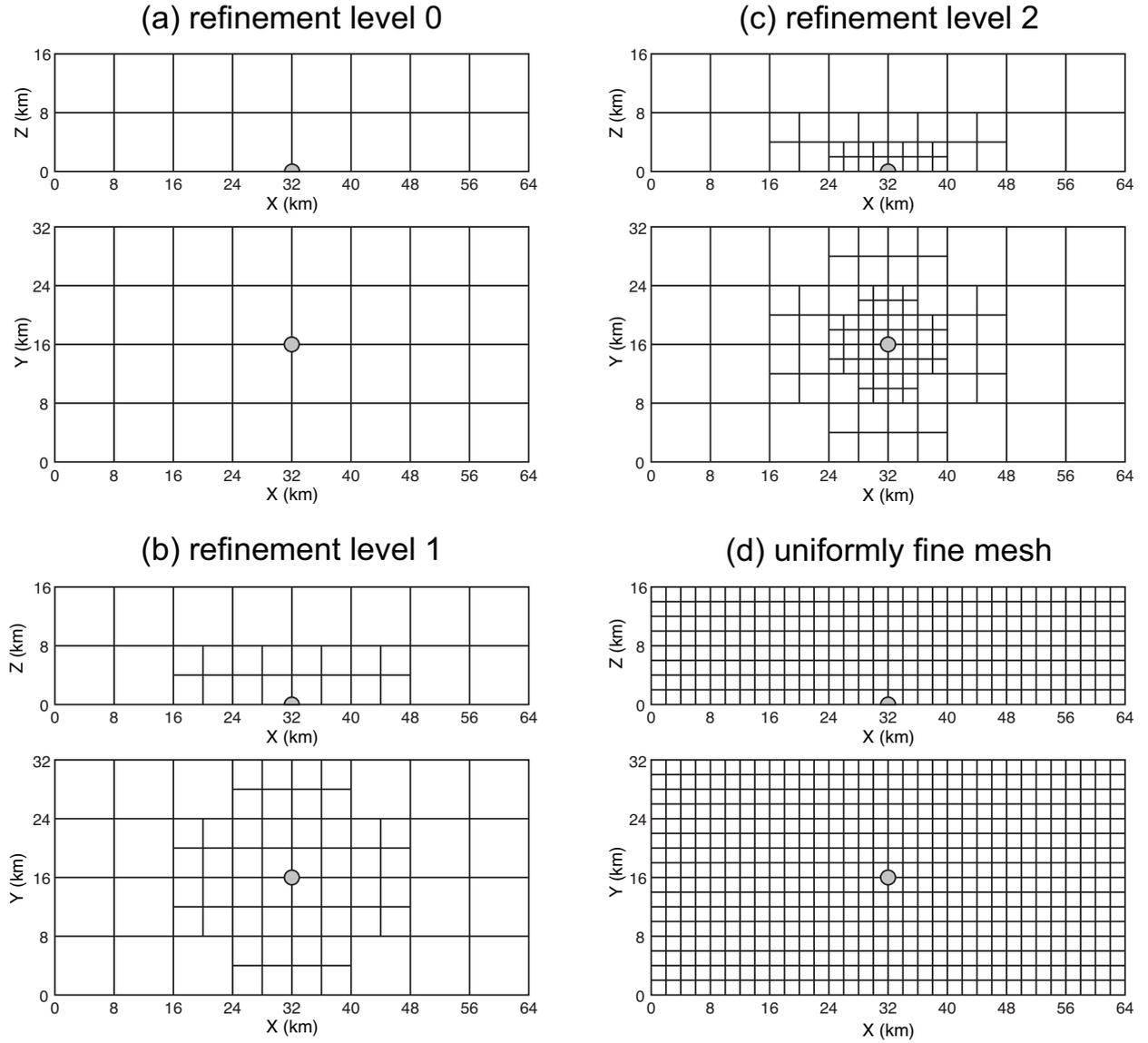}\\
  \caption{Computational grid used in the test of flow over the 3D hemisphere hill with refinement levels (a) 0, (b)1, (c) 2 and with (d) a uniformly fine mesh. The top and bottom panels in each figure describe the vertical and horizontal structures of the grid. Solid lines are the cube boundaries and each cube has 20 $\times$ 20 cells. The grid spacing in large, medium and small cubes is 400 m, 200 m, and 100 m, respectively. Shaded region describes the topography.}
  \label{Fig12}
\end{figure*} 

\begin{figure*}[t]
  \centering
  \noindent\includegraphics[width=39pc,angle=0]{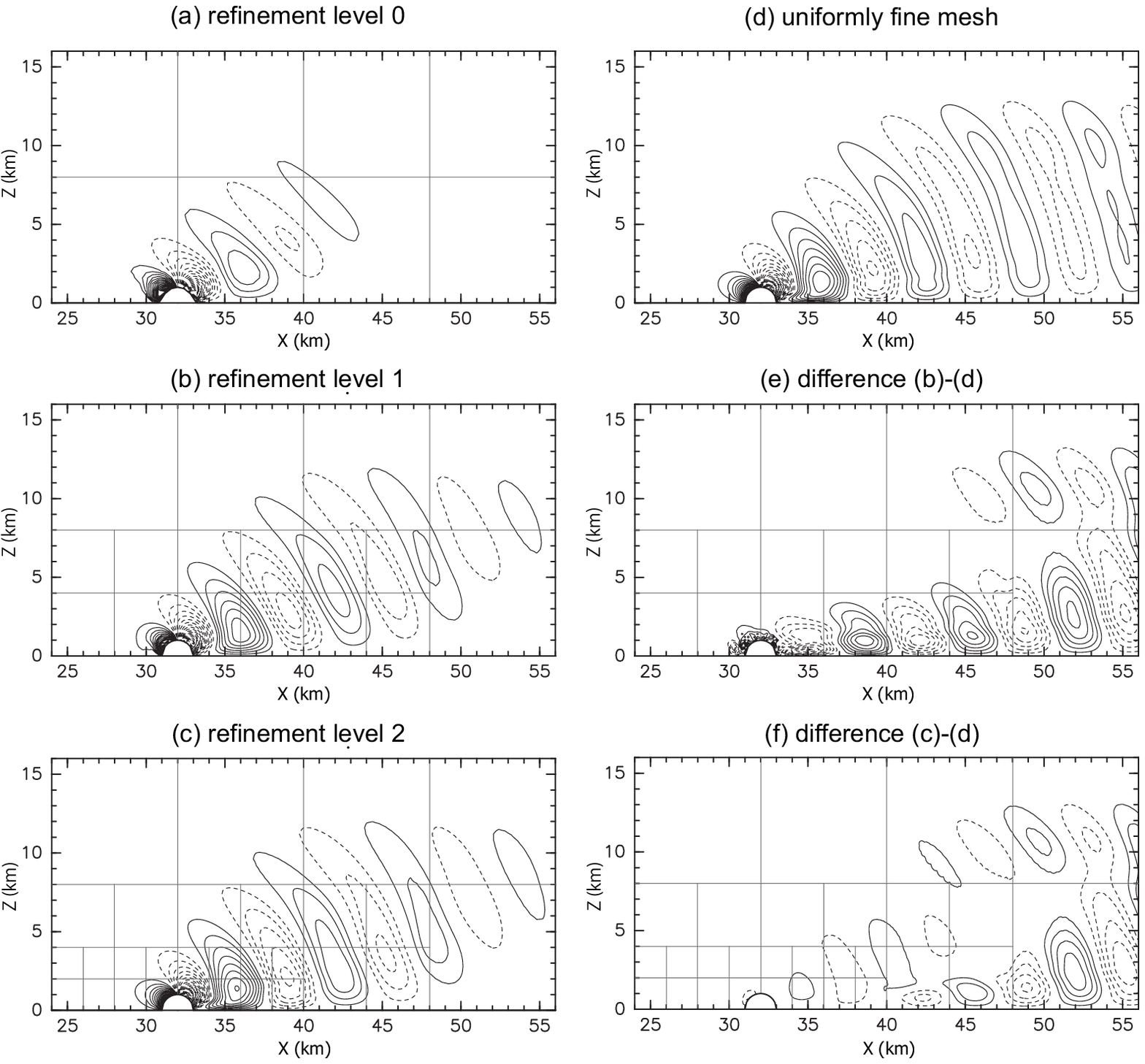}\\
  \caption{Vertical $x$-$z$ slices of the vertical velocity field through the center of the 3D hemisphere hill in the mesh refinement test case. Results are reproduced on the grid with refinement levels (a) 0, (b)1, (c) 2 and with (d) a uniformly fine mesh. Figure (e) shows the difference between the vertical velocity fields in (b) and (d), and (f) shows the difference between (c) and (d). Contour intervals are 0.25 m s$^{-1}$ in (a)--(d), and 0.1 m s$^{-1}$ in (e) and (f).
}\label{Fig13}
\end{figure*}

These results lead to the conclusion that the model reproduces the 3D structure of flow over the bell-shaped hill, including detailed features as predicted by the theory. In addition, it is shown that vertical combinations of cells in our model work properly over varying gradients.

\subsubsection{Flow over a 3D hemisphere-shaped hill} \label{ss:hemi}

The second test demonstrates the ability of the cut-cell model to simulate flow over very steep topography, using a hemisphere hill of radius $r$ = 1 km. As in the case with the bell-shaped hill, a uniform grid is used in this test with equal spacing of 200 m in all three directions. The maximum slope angle is resolved at 71$^{\circ}$ at the foot of the hill. Because the horizontal gradients of the terrain surface range from gentle around the top of the hill to very steep at the foot, all five directions of cell-merging described in section \ref{ss:cut-cells} are involved in this test.

Figure \ref{Fig11}(a) shows a vertical $x$-$z$ cross section of the simulated vertical velocity field. Due to the steepness of the topography, much larger amplitudes are calculated in the flow over the hemisphere hill compared to the result of Figure \ref{Fig10}(a) for the bell-shaped hill. Figures \ref{Fig11}(b) and (c) show the horizontal $x$-$y$ slices of the vertical velocity field over the hemisphere hill for heights of 1200 m and 2400 m, respectively. In the leeward side of the hill, damped oscillations towards downstream and U-shaped flow patterns widening with height are observed as in the case of the bell-shaped hill. Note that contour intervals in Figures \ref{Fig11}(b) and (c) are two and a half times larger than that in Figures \ref{Fig10}(b) and (c). 

\begin{figure*}[t]
  \centering
  \noindent\includegraphics[width=39pc,angle=0]{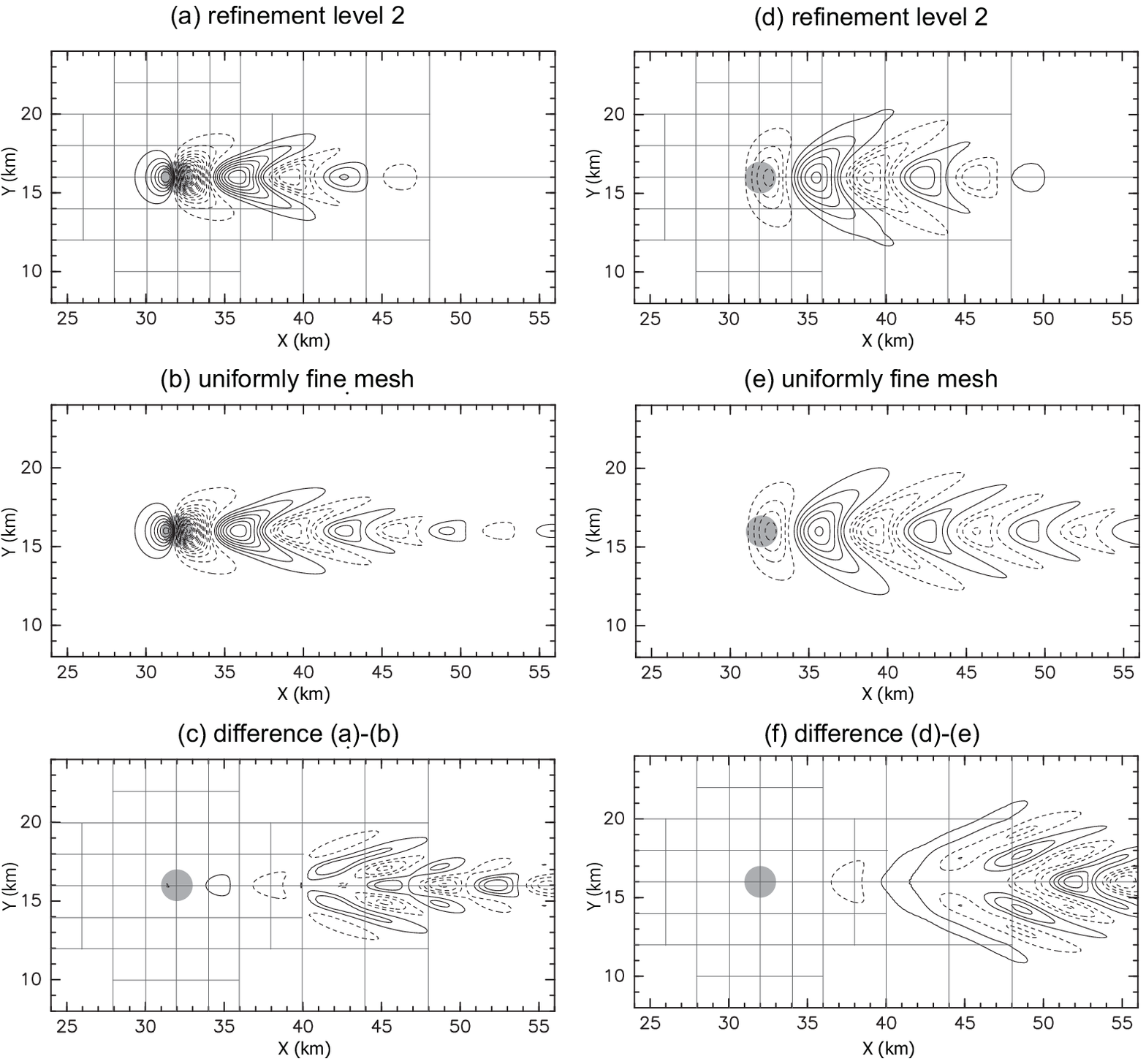}\\
  \caption{Horizontal $x$-$y$ slices of the vertical velocity field over the hemisphere hill at height $z$ = 1200 m reproduced on the grid with (a) the refinement level 2 and (b) a uniformly fine mesh. Figure (c) shows the difference between the vertical velocity in (a) and that in (b). Solid and dashed lines indicate positive and negative values, respectively. Contour intervals are 0.25 m s$^{-1}$ in (a) and (b), and 0.1 m s$^{-1}$ in (c). Gray lines in (a) and (c) describe the cube boundaries with the refinement level 2. Shaded region describe the topography. (d)--(f) same as with (a)--(c) except $z$ = 2400 m.}\label{Fig14}
\end{figure*}

\begin{figure}
  \centering
  \noindent\includegraphics[width=17pc,angle=0]{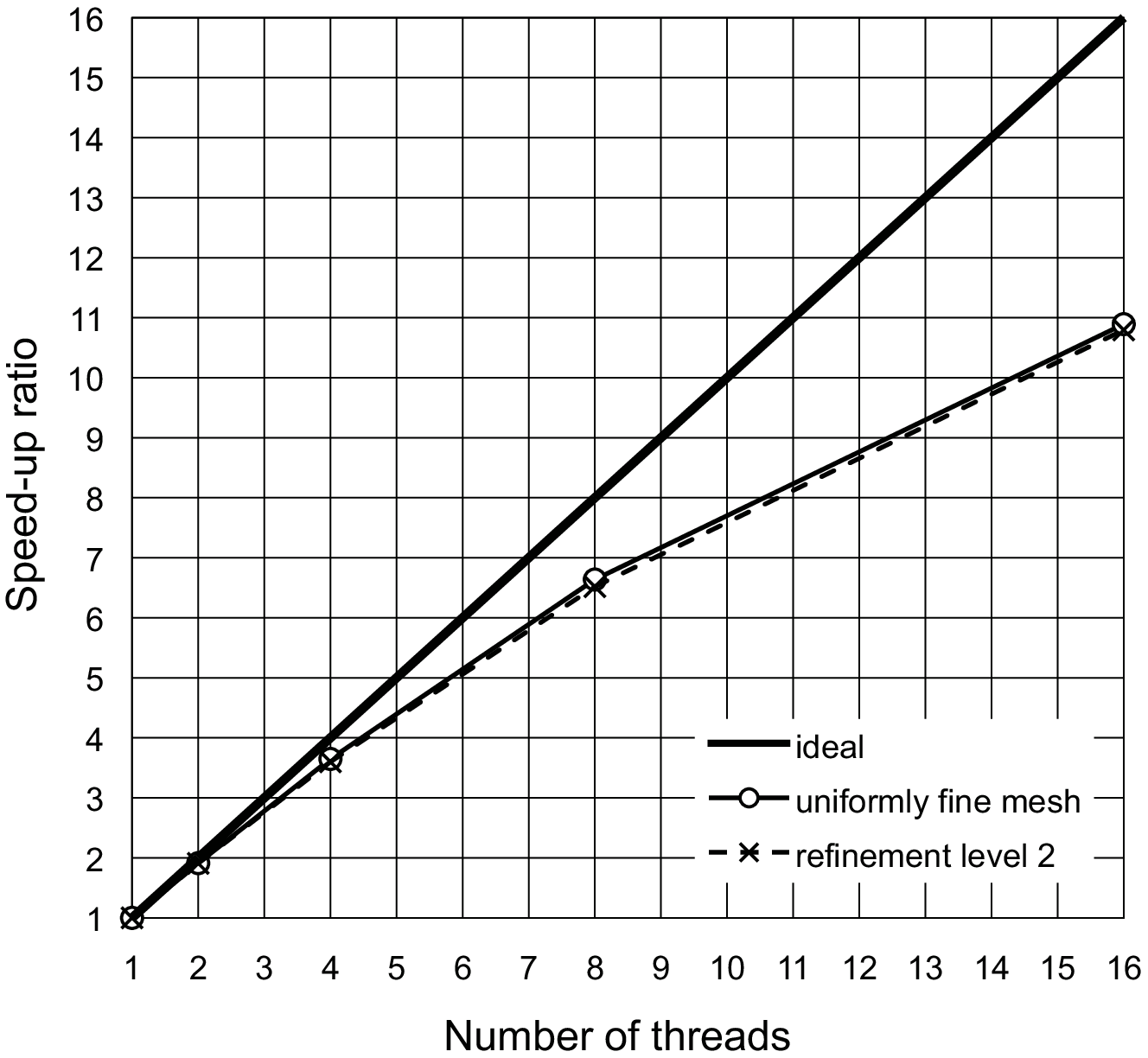}\\
  \caption{Comparison of speed-up ratio with increasing the number of threads.}
  \label{Fig15}
\end{figure} 

When the flow past an isolated sphere is symmetric about a horizontal plane through the centre, it can be equated to the flow past a hemisphere with the free-slip boundary condition at the surface \citep[][]{Baines1995}. Numerical results by \cite{Hanazaki1988} indicate that the flow past a sphere for the case with $Nr/U = 1$ is indeed symmetric about the horizontal plane. Referring to the isopycnic lines over the sphere calculated in \cite{Hanazaki1988} (see their Figures 3(c) and 11(f)), our model faithfully reproduces a smooth flow field close to the terrain surface, as illustrated in the result of the potential temperature on a vertical $x$-$z$ cross section through the center of the hill (Figure \ref{Fig11}(d)). \cite{LinEA1992} examined experimentally the lee wavelength in the flows past a sphere. The observed wavelength in the lee waves for the case with $Nr/U = 1$ was slightly less than $2\pi U/N \approx 6.3$ km, which is comparable to the apparent lee wavelength in our model result of Figure \ref{Fig11}(b) (approximately 7 km).  

Though the maximum slope angle of the hemisphere hill is almost as large as that of a very steep bell-shaped hill used in \cite{LockEA2012}, the reproduced flows show very different structures: vertically propagating mountain waves are reproduced in this test, whereas a simple flow that climbs over and goes down the hill is reproduced in \cite{LockEA2012}. These results can be explained in terms of the parameter $Na/U$ indicating the types of disturbances that would be observed in each case, according to linear theory \citep[e.g.,][]{Smith1977, SatomuraEA2003}. When $Na/U \gg 1$, the theory predicts that nearly hydrostatically balanced mountain waves will appear, while nonhydrostatic mountain waves should be observed when $Na/U \approx 1$. Evaluating $Na/U$ by $Nr/U = 1$ for the case of the hemisphere hill, the model results in Figure \ref{Fig11} show the generation of mountain waves as expected from the theory. Given the steep nature of the hemisphere hill, we may not expect the linear theory to give a good description of the flow. However, referring to the flow past a sphere when $Nr/U = 1$, which is classified by \cite{LinEA1992} as their flow regime ``WV'', our result appears to be comparable to the linear theory as discussed in \cite{Baines1995}. When $Na/U \ll 1$, buoyancy forces become less important and the motion approaches potential flow, such as observed in the case of \cite{LockEA2012}. 

We therefore conclude that our model successfully reproduces a smooth and physically reasonable flow over a hill which has very steep slopes. At the same time, it is proved that the model works properly even when vertically and horizontally merged cells coexist.

\subsection{Cut-cell model on a locally refined mesh}\label{s:refine}

Finally we demonstrate the performance of our cut-cell model using the block-structured mesh-refinement approach described in section \ref{ss:blocks}. In this test, four simulations of flow over the same 3D hemisphere-shaped hill used in section \ref{ss:hemi} are performed using grids with different refinement levels. 

The grids used in this test are shown in Figure \ref{Fig12}. Each figure shows the boundaries of cubes that comprise the grid. Top panels describe the cube boundaries on a vertical $x$-$z$ cross section through the center of the hill, and bottom panels describe those on a horizontal $x$-$y$ cross section at the sea level. Here every grid contains 20 $\times$ 20 cells in each cube. At the refinement level 0 (Figure \ref{Fig12}(a)), the domain is filled with large cubes of 400 m spacing in all three directions. At the refinement level 1 (Figure \ref{Fig12}(b)) and 2 (Figure \ref{Fig12}(c)), cubes are divided into finer cubes around the hill, where the medium- and small-sized cubes have spacings of 200 m and 100 m in all three directions, respectively. The results using the grids with each refinement level are compared to the result using a uniformly fine mesh with equal spacing of 100 m (Figure \ref{Fig12}(d)), which is considered as a reference. The total numbers of cubes are 64, 148, 232 and 4096 in Figures \ref{Fig12}(a), (b), (c) and (d), respectively. 

Figures \ref{Fig13}(a) to (d) respectively show vertical $x$-$z$ cross sections of the simulated vertical velocity fields using the grids with refinement levels 0, 1, 2 and using the uniformly fine mesh. Due to the lack of resolution, the amplitude of the vertical velocity using the grid with the refinement level 0 (Figure \ref{Fig13}(a)) decays faster with height compared to that using the uniformly fine mesh (Figure \ref{Fig13}(d)). In the leeward side of the hill, the propagating waves towards further downstream captured in Figure \ref{Fig13}(d) are missing in Figure \ref{Fig13}(a). The amplitude is successfully increased when the refined meshes are implemented around the hill (Figures \ref{Fig13}(b) and (c)). Though in the region of the coarse resolution the simulated amplitude is smaller than that in the reference solution of Figure \ref{Fig13}(d), both results using the grid with the refinement level 1 (Figure \ref{Fig13}(b)) and the refinement level 2 (Figure \ref{Fig13}(c)) successfully reproduce the wave propagation towards downstream. In addition, no visible distortions of mountain wave patterns are found in these results across the fine-coarse cube boundaries.

To compare the results of Figures \ref{Fig13}(b) and (c) in detail, the differences with respect to the reference solution of Figure \ref{Fig13}(d) are shown in Figures \ref{Fig13}(e) and (f), respectively. Note that the contour interval in Figures \ref{Fig13}(e) and (f) is two fifths of that in Figures \ref{Fig13}(a)--(d). In Figure \ref{Fig13}(e), some errors are found in both the regions at the refinement levels 0 and 1 where the resolution is coarser than that of the reference solution. These errors are successfully reduced when implementing the refined mesh at level 2 around the hill in Figure \ref{Fig13}(f). In particular, it is shown that refining the mesh around the hill not only improves the result in the refined region at level 2 but also improves the result in the coarser regions at levels 1 and 0. These results indicate that the benefit due to the higher resolution within the refined cubes exceeds the errors due to the resolution changes at fine-coarse cube boundaries.

Figure \ref{Fig14} displays horizontal $x$-$y$ slices of the simulated vertical velocity fields. Here Figures \ref{Fig14}(a) and (b) show the vertical velocity fields at $z$ = 1200 m using the grid with the refinement level 2 (Figure \ref{Fig12}(c)) and using the uniformly fine mesh (Figure \ref{Fig12}(d)), respectively. Figure \ref{Fig14}(c) shows the difference between the velocity values in Figures \ref{Fig14}(a) and (b), where the contour interval in Figure \ref{Fig14}(c) is two fifths of that used in Figures \ref{Fig14}(a) and (b). As in the vertical slices of the results, the model reproduces very accurate flow in the region at the refinement level 2 compared to the reference solution using the uniformly fine mesh. As the mesh resolution gets coarser towards the leeward side of the hill, the amplitude decays faster in Figure \ref{Fig14}(a) compared to that in the reference solution of Figure \ref{Fig14}(b), therefore the errors become larger as shown in Figure \ref{Fig14}(c). Because there are no visible distortion nor large errors found at the fine-coarse mesh boundaries, these errors are considered to be mainly due to the reduction of the mesh resolution. A similar behavior of the model is found at $z$ = 2400 m (Figures \ref{Fig14} (d) to (f)). Here Figures \ref{Fig14}(d) and (e) show the vertical velocity fields at $z$ = 2400 m using the grids of Figures \ref{Fig12}(c) and (d), respectively, and Figure \ref{Fig14}(f) shows the difference between them. These results confirm the model's ability to effectively increase the resolution around the hill without introducing large errors associated with the change of mesh resolution.

Finally the speed-up ratio is estimated in order to verify the scalability of the model both on a locally refined grid and on a uniform grid (Figure \ref{Fig15}). The parallelization of the model code is described in section \ref{ss:parallel}. Here the computations are performed on a Linux PC which has 16-core AMD Opteron CPUs.  The parallel efficiency is about 80\% and 70\% at 8 and 16 processors, respectively, in both tests using the grid with the refinement level 2 (Figure \ref{Fig12}(c)) and using the uniformly fine mesh (Figure \ref{Fig12}(d)). This result indicates that the mesh refinement and the subcycling time integration do not affect the parallel efficiency of the overall computation in our model.

These results lead to the conclusion that the model reproduces a flow result on a locally refined grid with both sufficient accuracy and computational efficiency, thereby demonstrating the advantage of the model's mesh-refinement technique in combination with the cut-cells for high-resolution atmospheric simulations over 3D topography. 

\section{Conclusions} \label{conclusion}

To achieve high-resolution and highly precise simulations over 3D topography, a new 3D nonhydrostatic atmospheric model was developed using the Cartesian cut-cell method. To avoid severe restriction on time steps while maintaining the sharp representation of terrain surface, merged cells are used in both the horizontal and vertical directions. In addition, a block-structured mesh-refinement approach was implemented to allow high resolution at the boundary layers with computational efficiency. 

The accuracy of the model was confirmed by simulating flow over a 3D bell-shaped hill and comparing the result with the analytical solution given by the linear theory. The ability of the model to simulate flow over topography including a very steep slope was demonstrated using a hemisphere-shaped hill. The model successfully reproduced smooth and physically reasonable flow including detailed features which agree with numerical and experimental results. Finally, a test of the mesh-refinement approach in combination with the cut-cell representation of topography was performed using the same hemisphere-shaped hill. The model reproduced smooth mountain waves propagating over varying grid resolution without introducing large errors at fine-coarse cube boundaries. In addition, the model showed a good parallel efficiency on a locally refined grid around the hill.

A remaining issue with the use of Cartesian cut-cell grids is the incorporation of physical process models and other components, such as a land surface model. As many of the physical parameterization schemes are based on a single column model, the couplings with the block-structured grid and with horizontally merged cells may not be straightforward. Further work is intended to solve these issues with the aim of simulating real-world atmospheric problems. 

\begin{acknowledgements} 
We acknowledge the contributions of the second author, Prof. Takehiko Satomura, who sadly passed away during the preparation of the manuscript. Without his enthusiasm, encouragement and our discussions over many years this work would have never been completed. The first author would like to thank Dr Colin J. Cotter for his kind support. We also thank Mr Ed White and two anonymous reviewers for their helpful comments.

This work was supported by the JSPS Postdoctoral Fellowship for Research Abroad. Part of this work was supported by a Grant-in-Aid for JSPS Bilateral Joint Research Project. Some figures were drawn by the GFD-DENNOU Library.
\end{acknowledgements}

\appendix

\section{Estimation of Global Accuracy}
In this section, we numerically demonstrate that our cut-cell method is globally second-order accurate for a flow over 3D terrain. Here we performed four simulations of flow over a 3D cosine-shaped hill using a uniform grid with different grid intervals: $\Delta$ = 400 m, 250 m, 200 m, and 100 m. The same grid interval is used in $x$, $y$ and $z$ directions for each test, and the same atmospheric and boundary conditions as used in section \ref{s:result} are imposed. The shape of the hill is described by
\begin{eqnarray}
	h(x,y) &=& \left\{ \begin{array} {lcl} \dfrac{h_{\mathrm{m}}}{2} \left[1+\cos\left(\dfrac{\pi r}{a}\right)\right] & \mathrm{if}  & r < a \\  \noalign{\vskip 1.5ex}
	0 & \mathrm{if} & r \geqq a, \\ \end{array} \right. 
\end{eqnarray}
where 
\begin{eqnarray}
r = \sqrt{(x-x_{c})^{2}+(y-y_{c})^{2}}
\end{eqnarray}
denotes the distance between a position ($x$, $y$) and the center of the lower boundary ($x_{c}$, $y_{c}$): in this test, $x_{c}$ = 32 km and $y_{c}$ = 16 km. The maximum height $h_{\mathrm{m}}$ and the half-width $a$ are set to 400 m and 2 km, respectively.

\begin{figure}[t]
  \centering
  \noindent\includegraphics[width=15pc,angle=0]{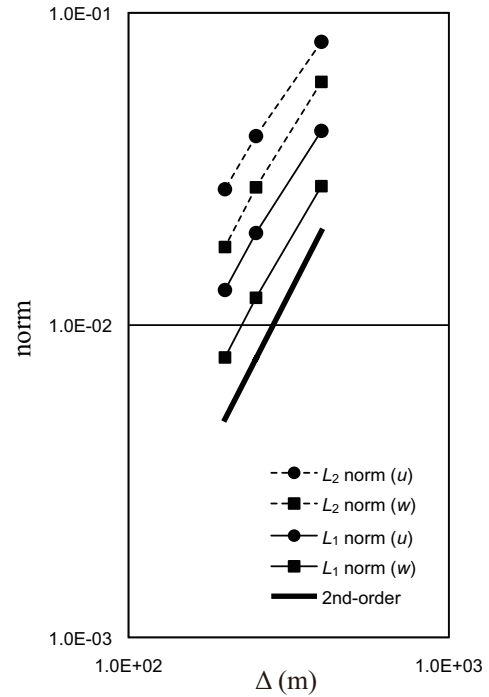}\\
  \caption{Variation of global error in $u$ and $w$ with grid intervals for flow over the cosine-shaped hill. Thin solid and dashed lines indicate $L_{1}$ and $L_{2}$ norms of the errors, respectively, and thick solid line corresponds to second-order accurate convergence.}
  \label{FigA1}
\end{figure}  

Following \cite{Yamazaki+Satomura2010}, both the $L_{1}$ and $L_{2}$ norms of the global errors are computed as,
\begin{flalign}
	\varepsilon_{L_{1}} &= \frac{1}{N_{x}N_{y}N_{z}}\sum_{I=1}^{N_{x}N_{y}N_{z}}|\psi_{I}^{numerical}-\psi_{I}^{exact}| \hspace{10pt}, \\
	\varepsilon_{L_{2}} &= \left(\frac{1}{N_{x}N_{y}N_{z}}\sum_{I=1}^{N_{x}N_{y}N_{z}}(\psi_{I}^{numerical}-\psi_{I}^{exact})^{2}\right)^{1/2},  
\end{flalign}
respectively, by assuming the solution with the finest solution of $\Delta$ = 100 m as the exact solution. Here $N_{x}$, $N_{y}$ and $N_{z}$ indicate the grid numbers above the top of the hill in $x$, $y$ and $z$ directions, respectively. Figure \ref{FigA1} shows a log-log plot of the errors in velocity components $u$ and $w$ versus grid intervals. Also shown is a line with a slope of 2 which corresponds to the second-order accurate convergence. The plot clearly shows that the global error in our computed solution decreases in a manner consistent with a second-order accurate scheme. This test therefore shows that our method produces results which are consistent with a method of global second-order accuracy.
 

\bibliographystyle{plainnat}
\bibliography{Yamazaki_references}
\end{document}